\documentclass[11pt,twoside]{article}
\usepackage{latexsym,amssymb,amsmath}
\RequirePackage{graphics,epsfig,psfrag}%,picins}
\usepackage{graphbox}
\usepackage[latin1]{inputenc}
\usepackage{color}
\usepackage{colordvi,fancyhdr}
\usepackage{pdfpages}
\definecolor{Gris}{cmyk}{0.1,0.1,0.1,.75}
\usepackage[colorlinks=true,citecolor=Gris,linkcolor=Gris,filecolor=Gris,urlcolor=Gris]{hyperref}
\def\abs#1{\left \vert #1 \right \vert}

\def\RR{{\bf R}} %reelle Zahlen
\def\ZZ{{\bf Z}} %reelle Zahlen
\def\CC{{\bf C}} %reelle Zahlen
\def\LL{{\mathcal{L}}} %reelle Zahlen
\def\Mod#1{\,(\hbox{\rm mod}\,#1)}

\def\H{{\mathrm{H}}}
\def\TT{{\rm T}}

\def\phi{\varphi}
\def\eps{\varepsilon}

\def\bb{{\mathfrak{b}}}

\def\lk{\mathrm{lk}\,}
\def\pn{\medskip\par\noindent}
\def\bi{\vspace{-6pt}\begin{itemize}\itemsep -2pt plus 1pt minus 1pt}
\def\ei{\end{itemize}\vspace{-4pt}}
\def\bn{\vspace{-6pt}\begin{enumerate}\itemsep -2pt plus 1pt minus 1pt}
\def\en{\end{enumerate}\vspace{-4pt}}
\newcommand{\Pf}{{\em Proof}. }

\newcommand{\EPf}{\hbox{}\hfill$\Box$\vspace{.5cm}}
\def\[#1\]{\begin{eqnarray}#1\end{eqnarray}}
\def\$#1\${\begin{eqnarray*}#1\end{eqnarray*}}
\def\pent#1#2{\lfloor\frac{#1}{#2}\rfloor}

\def\sign#1{{\rm sign}\,\bigl( #1 \bigr)}
\def\Sum{\mathop{\sum}\limits}

\def\abs#1{\left \vert #1 \right \vert}

\def\frac#1#2{{\textstyle{{#1} \overwithdelims.. {#2}}}}
\def\Frac#1#2{{\displaystyle{{#1} \overwithdelims.. {#2}}}}
\def\@begintheorem#1#2#3{\par\addvspace{8pt plus3pt minus2pt}%
              \noindent{\csname#1headfont\endcsname#1\ \ignorespaces#3 #2.}%
              \csname#1font\endcsname\hskip6pt\ignorespaces}
\def\@endtheorem{\par\addvspace{8pt plus3pt minus2pt}\@endparenv}

\catcode`\@=12

\newtheorem{theorem}{Theorem}[section]
\newtheorem{thm*}{Theorem}
\newtheorem{thm}[theorem]{Theorem}
\newtheorem{corollary}[theorem]{Corollary}
\newtheorem{lemma}[theorem]{Lemma}
\newtheorem{proposition}[theorem]{Proposition}
\newtheorem{definition}[theorem]{Definition}
\newtheorem{remark}[theorem]{Remark}

\newtheorem{example}[theorem]{Example}

\def\[#1\]{\begin{align}#1\end{align}}
\begin{document}
\date{\today}
%\pagestyle{myheadings}
%\markboth{E. Brugall\'e, P. -V. Koseleff, D. Pecker}
%{The lexicographic degree of the first two-bridge knots}
\title{The lexicographic degree of the first two-bridge knots}
%{The lexicographic degree \\ of the first two-bridge knots}
\author{Erwan Brugall\'e, Pierre-Vincent Koseleff, Daniel Pecker}
\maketitle
%    General info
%pvk 11/09 \linenumbers
%pvk 11/09 \modulolinenumbers[5]
%%%%%%%%%%%%%%%%%%%%%%%%%%%%%%%%%%%%%%%%%%%%%%%%%%%%%%%%%%%%%%%%%%%%%%%%%%%
\begin{abstract}
We study the degree of polynomial representations of knots.
We give %here
the lexicographic degree of all two-bridge knots
with 11 or fewer crossings.
First, we estimate the total degree of a
lexicographic parametrisation of such a knot.
This allows us to transform this problem into a study of real
algebraic trigonal plane curves,
%eb
and in particular to
%. We
use the braid theoretical method developed by Orevkov.
\\[10pt]
\\[10pt]
\end{abstract}
\begin{center}
\parbox{12cm}{\small
\tableofcontents }
\end{center}
\section{Introduction}\label{sec:intro}
\def\imagetop#1{\vtop{\null\hbox{#1}}}
A polynomial parametrisation of
 a knot $K$ in $ {\bf S^3}$ is a polynomial
map $\gamma:\RR\to\RR^3$ whose closure of the image in ${\bf S^3}$ is
 isotopic to $K$.
Every knot admits a polynomial parametrisation, see \cite{Sh,Va}.
In this paper we are interested in determining the
\emph{lexicographic degree} of a knot $K\subset {\bf S^3}$, i.e. the minimal
degree for the lexicographic order of a polynomial
parametrisation of $K$.

\begin{figure}[!ht]
\begin{center}
\begin{tabular}{ccc}
{\scalebox{.8}{\includegraphics{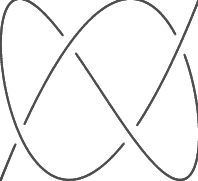}}}&&
{\scalebox{.8}{\includegraphics{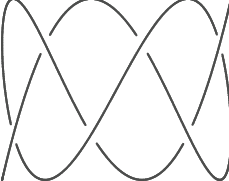}}}\\
$\deg 4_1 =(3,5,7)$&&$\deg 5_1 =(3,7,8)$
\end{tabular}
\end{center}
\caption{\small  Trigonal polynomial
  diagrams of the figure-eight knot $4_1$ and the torus knot $5_1$}
\label{fig:h357}
\end{figure}

\begin{figure}[!ht]
\begin{tabular}{ccc}
 {\scalebox{.5}{\includegraphics{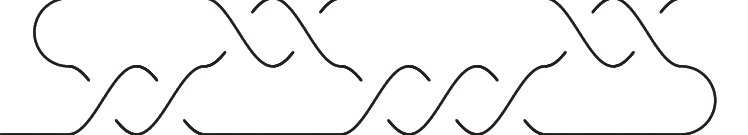}}} &\quad &
 {\scalebox{.5}{\includegraphics{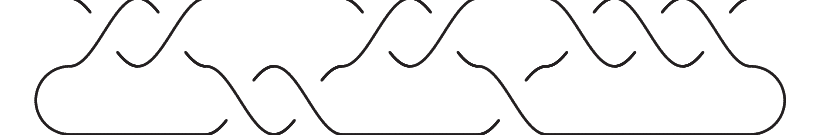}}}\\
{\small A 9-crossing diagram of degree $\geq (3,13,14)$} & &
{\small A 10-crossing diagram of degree $(3,11,16)$}
\end{tabular}
\caption{\small  Two diagrams of $9_{15}$}
\label{fig:9_15}
\end{figure}

The unknot has
lexicographic degree $(-\infty,-\infty,1)$, and it is easy to see that
%one sees easily that
the  lexicographic degree of any other knot is
$(a,b,c)$ with $3\le a<b<c$.
\emph{Two-bridge} knots are precisely those with lexicographic degree
$(3,b,c)$, see \cite{KP4};
they have a $xy$-projection which is a trigonal curve.
See
Figure \ref{fig:h357} for two examples of trigonal polynomial parametrisations
of a long knot.
\pn
Two-bridge knots are an important family of knots. The first 26 knots (except $8_5$) are two-bridge knots. Moreover these knots are classified by their Schubert fractions, which can be easily computed from any trigonal projection, see Section \ref{sec:2-bridge}.
%The authors have already obtained preliminary results in
%\cite{BKP1}, focusing on arithmetic properties of continued fractions and
%in \cite{BKP2}, focusing on topological properties of trigonal plane curves.
\pn
One might expect that the lexicographic degree of a knot $K$ is obtained for a minimal-crossing diagram of
this knot. This is not true.
The diagram on the left of Figure \ref{fig:9_15} is a minimal crossing diagram
%an alternating (hence minimal-crossing) diagram
of the knot  $ 9_{15}$.
On the right of the figure is a
10-crossing diagram of smaller degree of the same knot.
This is why it is necessary to consider all the diagrams of two-bridge knots.
The enumeration of all possible diagrams of a given two-bridge knot can be
efficiently done using Conway's notation.
\pn
In this paper, we show:
\pn
{\bf Theorem.}
{\em The lexicographic degree of all 186 two-bridge knots with crossing
number at most $N\le 11$ is $(3,b,3N-b)$, where the values of $b$ are listed in
Table \ref{table:all}, p. \pageref{table:all}.}
\pn
We prove this result in two steps.
\pn
{\bf Proposition \ref{prop:N<=11}}
{\em The lexicographic degree $(3,b,c)$ of a knot with crossing number
$N\le 11$ satisfies $b+c=3N$.}
\pn
Proposition \ref{prop:N<=11} also holds for all $N$ when $b \leq N+3$ or
$b = \pent{3N-1}2$.
We prove in Theorem \ref{thm:b+c} that $b+c\ge 3N$ for any
polynomial parametrisation of degree $(3,b,c)$ of a knot with crossing
number $N$. Furthermore,
every two-bridge knot of crossing number $N$ admits a parametrisation of
degree $(3,b,c)$ with $b+c=3N$, see \cite{KP4}.
We do not know if Proposition \ref{prop:N<=11} holds for all crossing numbers $N\geq 12$.
\pn
Proposition \ref{prop:N<=11} allows us to reduce the determination
of the lexicographic degree of a two-bridge knot to the study of plane
curves.
For knots with 11 crossings or fewer, it is enough to determine the
smallest integer $b$ such that a plane projection
admits a polynomial parametrisation of degree $(3,b)$. This reduction
to plane curves enlarges the set of tools at our disposal; in
particular we make an important use of Orevkov's braid theoretical
approach in the study of pseudoholomorphic curves.
\pn
Hence the second step in the proof of our theorem
is to focus on parametrisations of plane projections.
We introduce the $\TT$-reduction in Section \ref{sec:upper},
that corresponds to the projection of the Lagrange isotopy on
trigonal diagrams. The $\TT$-reduction allows us to remove a triangle of crossings
from a diagram,
and therefore to obtain an upper bound for degrees we are looking for.
On the other hand, we introduce the $\TT$-augmentation in Section \ref{sec:augmentation} that allows us to add
a triangle of crossings to a given diagram $D$. From a polynomial parametrisation corresponding to $D$
we deduce a parametrisation for the new diagram.
\pn 
We propose an algorithm to find the lexicographic
degrees of the first 186 two-bridge knots with 11 crossings or fewer.
As a byproduct of our computations, we also exhibit in Table \ref{table:noalt}
the 16 two-bridge knots with 11 crossings or fewer for which the lexicographic degree
is smaller than the degree of their minimal-crossing diagrams.
\pn
The paper is organised as follows.
In Section \ref{sec:2-bridge} we recall Conway's notation for
trigonal diagrams of two-bridge knots.
Then we prove the inequality $b+c \geq 3N$ in Section \ref{sec:lower}
and deduce Proposition \ref{prop:N<=11}. % below.
In Section \ref{sec:plane}, we consider
plane trigonal curves and we first obtain a lower bound
for the lexicographic degree of a trigonal polynomial embedding
in Proposition \ref{prop:Bezout}.
We obtain another bound for pseudoholomorphic curves
and therefore for polynomial embeddings in Proposition \ref{prop:even}.
In Section \ref{sec:lex}, we obtain
the lexicographic degrees of the first 186 two-bridge knots
with 11 crossings or fewer.

\section{A lower bound for the total degree
of two-bridge knots}\label{sec:b+c}

\subsection{Trigonal diagrams  of two-bridge knots}\label{sec:2-bridge}
\begin{figure}[!ht]
\begin{center}
{\scalebox{.75}{\includegraphics{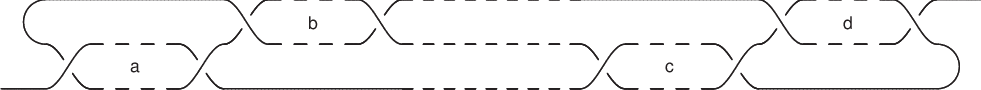}}}\\[30pt]
{\scalebox{.75}{\includegraphics{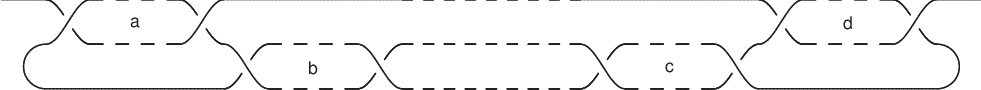}}}
\end{center}
\caption{\small  Conway's  form for two-bridge knots (or links)}
\label{fig:conways3}
\end{figure}

\begin{figure}[!ht]
\begin{center}
\begin{tabular}{ccc}
{\scalebox{.9}{\includegraphics{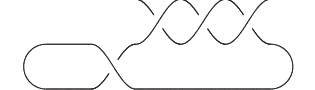}}}&\quad&
{\scalebox{.9}{\includegraphics{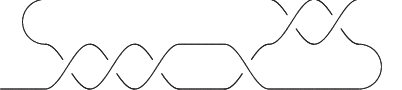}}}\\
$C(0,1,3)$&& $C(3,0,-1,-2)$
\end{tabular}
\end{center}
\caption{\small  Examples of trigonal diagrams}
\label{fig:conways4}
\end{figure}

\begin{figure}[!ht]
 \centerline{
 \psfrag{a}{$m-1$}\psfrag{b}{$1-n$}
 {\scalebox{.70}{\includegraphics{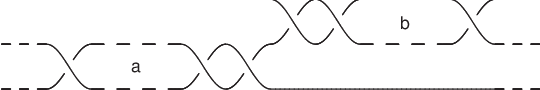}}} \quad
\quad \quad
 \psfrag{a}{$m-1$}\psfrag{b}{$n-1$}
 {\scalebox{.70}{\includegraphics{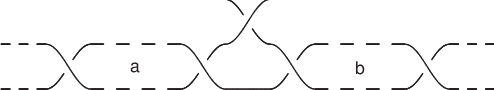}}}
}
 \caption{\small  Lagrange isotopy: $C(u,m,-n,-v) \to C(u,m-1,1,n-1,v)$}
 \label{fig:lagrange}
 \end{figure}

A two-bridge knot  admits a diagram in \emph{Conway's  open form}
(or trigonal form). This diagram, denoted by
$C(m_1, m_2, \ldots, m_k)$  where $m_i  \in \ZZ$,
is explained by Figure \ref{fig:conways3} (see \cite{Co}, \cite[p.~187]{Mu}).
\psfrag{a}{\small $m_1$}\psfrag{b}{\small $m_2$}
\psfrag{c}{\small $m_{k-1}$}\psfrag{d}{\small $m_{k}$}
The number of twists is denoted by the integer $\abs{m_i}$, and the sign of
$m_i$ is defined as follows: if $i$ is odd, then the right twist is positive,
if $i$ is even, then the right twist is negative.
In Figure \ref{fig:conways3}  the integers $m_i$ are all positive.
 Figure \ref{fig:conways4} shows the examples $C(0,1,3)$, $C(3,0,-1,-2).$
\pn
The two-bridge knots (or links) are classified by their Schubert fractions
$$
\Frac {\alpha}{\beta} =
m_1 + \Frac{1} {m_2 + \Frac{1} {\cdots +\Frac 1{m_k}}}=
[ m_1, \ldots, m_k], \quad \alpha \geq 0, \ (\alpha,\beta)=1.
$$
Given $[ m_1, \ldots, m_k]=\Frac {\alpha}{\beta} $ and
$[ m'_1, \ldots, m'_l]=\Frac {\alpha'}{\beta'}$, the diagrams
 $C(m_1, m_2, \ldots, m_k)$ and $C(m'_1, m'_2, \ldots, m'_l)$
correspond to isotopic knots (or links) if and only if
$\alpha = \alpha' $ and $ \beta' \equiv \beta ^{\pm 1} \Mod{\alpha}$,
see \cite[Theorem 9.3.3]{Mu}.
\pn
Every positive fraction ${\alpha}/{\beta}$ admits a continued fraction expansion
$[m_1, \ldots, m_k]$ where all the $m_i$ are positive.
Therefore every two-bridge knot $K$
admits a diagram in \emph{Conway's normal form}, that is an alternating
diagram of the form $ C(m_1, m_2, \ldots m_k)$,
where the  $m_i$ are all positive or all negative.
\pn
It is classical that one can transform any trigonal diagram
of a two-bridge knot into  Conway's normal form using the
Lagrange isotopies,  see \cite[p.~204]{Cr}.
\begin{definition}
Let $C(u,m,-n,-v)$ be a trigonal diagram, where $m,n$ are integers,
and $u,v$  are (possibly empty) sequences of
integers, see Figure \ref{fig:lagrange}. The Lagrange isotopy on $D$ is
\[
C( u,m,-n,-v) \rightarrow  C(u,m-\eps,\eps,n-\eps,v), \ \eps=\pm 1, \label{form:lagrange}
\]
\end{definition}
\pn
If $D=C(m_1, \ldots, m_k)$ is not in Conway's normal form, then it may happen that
$m_1=0$ or $m_k=0$. In this case, the diagram $D'=C(m_3, \ldots, m_k)$
or $D'=C(m_1, \ldots, m_{k-2})$ respectively, is the {\em reduced diagram} of $D$.
Since the diagram $C(m_1, \ldots, m_{i},0,0, m_{i+1}, m_{i+2}, \ldots, m_k)$ is identical
to $C(m_1, \ldots,  m_k)$, we can assume that if $m_i = 0$ then $m_{i-1} m_{i+1} \not =0$.
\pn
Given a finite integer sequence
$(m_1, \ldots, m_k)$, we  say that there is a {\em sign change}
between $m_i$ and $m_{i+1}$ if $m_{i} m_{i+1}<0$ or if $m_i=0$ and
$m_{i-1} m_{i+1} <0$.
\begin{proposition}\label{prop:Nlessthan}
Let $C (m_1, \ldots , m_k) $ be a diagram of a knot with crossing number
$N$. Let $N_0 = \Sum_{i=1}^k \abs{m_{i}}$ be the number of crossings,
and $\sigma$ be the number of sign changes in the sequence
$(m_1, \ldots , m_k)$. Then we have
$$
N \le  N_0  -\sigma.
$$
\end{proposition}
\Pf
If $\sigma=0$, then the inequality  means that the crossing number of a knot
is not greater than the number of crossings
of  a diagram of this knot.
Consequently, we can suppose $\sigma \ge 1.$
Let us prove the result by induction on
$N_0 = \Sum_{i=1}^k \abs{m_i} .$
We have to consider two cases.
\pn
First, let us suppose that the diagram  is of the form
$ C ( u, m, -n , -v )$, $m,n >0$
Then by a Lagrange isotopy we see that
$   C(u, m-1, 1, n-1 ,v) $
is another diagram of $K.$
In this new diagram, the number of crossings and the number of sign changes
are both  diminished by $1$.
Therefore we obtain by induction:
$$
 N  \le   (N_0 -1)  - (\sigma-1) = N_0 -\sigma  .
$$
Next, let us consider a diagram of the form
$ C ( u, m, 0 , -n ,v) ,  \  mn >0 .$
In this case we consider the new diagram
$ C( u, m-n, v).$
If $\sigma'$ is the number of sign changes of this new diagram,
then a
case by case
inspection shows that $ \sigma' \ge \sigma -2.$
As the number of crossings is diminished by at least $2$,
we obtain by induction:
$$
 N \le   ( N_0 -2)  - ( \sigma-2) = N_0 -\sigma    ,
$$
which concludes the proof.
\EPf
\pn
The proof of Proposition \ref{prop:Nlessthan} also implies the
following lemma.
\begin{lemma}\label{rem:N<N0}
In the notation of Proposition \ref{prop:Nlessthan}, we have:
\bn
\item If $\sigma=0$, then $N<N_0$ if and only if $m_1 \cdot m_k=0$.
\item If $\sigma=1$, then we have $N<N_0-1$ if and only if one of
the following situations occurs:
\bi
\item[--] $m_1=0$ or $m_k=0$,
\item[--] there exists $i$ such that $m_i=0$ and $m_{i-1}m_{i+1}<0$,
\item[--] $\abs{m_1}=1$ and $m_1 m_2 <0$ or $\abs{m_k}=1$ and $m_{k-1} m_k<0$.
\ei
\en
\end{lemma}
Let $D$ be a  long knot diagram, and
$\gamma:\RR\to\RR^3$ be  a parametrisation of $D$  whose crossing points
corresponds to the parameters $t_1<  \cdots < t_{2m}$.
Recall that the {\em Gauss sequence} of $D$  is the sequence
$g_1,\ldots,g_{2m}$
where $g_i=1$ if $t_i$ corresponds to an overpass, and $g_i=-1$
otherwise.
\begin{proposition}\label{prop:s}
Let  $C(m_1, \ldots, m_k)$, $m_i \not =0$, be a trigonal diagram of a
knot $K$, and
 $N_0 =\Sum \abs{m_i}$.
Let $s$ be the number of sign changes in the Gauss sequence of the diagram,
$\sigma$ be the number of sign changes  in the sequence
$(m_1, \ldots, m_k)$, and $\sigma_2$ be the number of
consecutive sign changes in the sequence $(m_1, \ldots, m_k)$.
Then, we have $$s = 2 N_0 - 3 \sigma + 2 \sigma_2 -1.$$
\end{proposition}
\Pf
We proceed
by induction on $(\sigma_2,\sigma)$.
If $\sigma=0$ then $\sigma_2 = 0$ and the diagram of $K$ is
alternating. In this case we have $s=2 \Sum \abs{m_i}-1 = 2 N_0 -1$.
\pn
If $\sigma_2=0$, we may assume that
$m_1>0$.
Let $j$ be the first index just that $m_i<0$. Then $j=k$ or
$m_{j+1}< 0$,
because $\sigma_2=0$. Let us consider the knot $K'$ defined by
$K' = C(m_1, \ldots, m_{j-1}, - m_{j}, -m_{j+1}, \ldots, -m_k)$.
We see that the number of sign changes in the Conway sequence of $K'$ is
$\sigma' = \sigma -1$, and that we still have $\sigma'_2 =0$.
By induction we get
$s'+3{\sigma'} = 2 \Sum \abs{m_i}-1$.
Since we have
$s'=s+3$, this completes the proof
when $\sigma_2=0$.
\pn
Now, let us suppose that $\sigma_2>0$ and consider
the first index $j$ such that $m_{j-1}\ m_{j}<0$ and $m_{j} \ m_{j+1} <0$.
Consider $K'$ defined by
$K' = C(m_1, \ldots, m_{j-1}, - m_{j}, -m_{j+1}, \ldots, -m_k)$.
We see that the number of sign changes in the Conway sequence of $K'$ is
$\sigma' = \sigma -1$ and also $\sigma'_2 =\sigma_2-1$.
By induction we get $s'+3\sigma' -2 \sigma'_2= 2 \Sum \abs{m_i}-1$.
Since we have $s'=s+1$, this
concludes the proof.
\EPf
\subsection{Total degree of two-bridge knots}\label{sec:lower}
The next theorem provides a
lower bound on the total degree of every trigonal knot diagram.
It generalises \cite[Theorem 4.3]{BKP2},
which proves that the lexicographic degree of
a knot of crossing number $N$
%$K$
is at least $(3,N+1,2N-1)$.
\begin{thm}\label{thm:b+c}
Let $\gamma :  \RR \rightarrow \RR^3 $ be a  polynomial parametrisation
of degree $(3, b,c)$ of a  knot of crossing number $N$. Then we have
$$b+c\ge  3N.$$
\end{thm}
\Pf We shall denote our polynomial knot $ \gamma (t) = ( x(t), y(t), z(t)) .$
 Without loss of generality, we may  assume that $b$ is not divisible
 by $3$.
 Let
 $C( m_1, m_2, \ldots, m_k)$ be the corresponding $xy$-diagram.
 To simplify the exposition, we shall first
suppose that $ m_i \ne 0$ for $ i=2, \ldots , k-1.$

By the genus formula,
the plane curve $C$ parametrised by
$C(t)= (x(t),y(t))$ has exactly $b-1$ nodes in $\CC^2$.
Let $N_0= \sum_{i=1} ^k  | m_i| $ be the number of real crossings of
$C$
(i.e. real nodes of $C$ which are the intersection of two real
branches of $C$),
and let $ \delta= b-1- N_0$ be the number of other nodes of $C.$

The real crossings are ordered by increasing abscissae. A real crossing is called \emph{special} if its Conway sign (for the trigonal
diagram) is  different from the Conway sign of the
preceding crossing.
\begin{figure}[!th]
\begin{center}
\begin{tabular}{ccc}
{\scalebox{.5}{\includegraphics{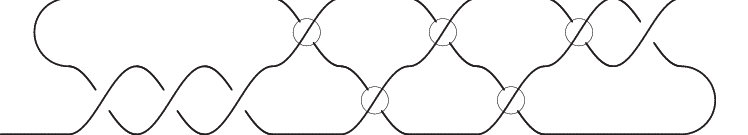}}}&&
{\scalebox{.5}{\includegraphics{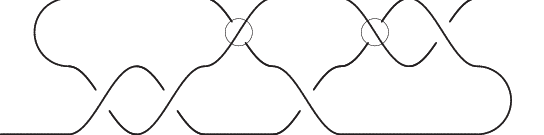}}}
\end{tabular}
\end{center}
\caption{\small  Special crossings of $C(3,-1,1,-1,1,-2)$ and $C(2,-1,-1,2)$}
\label{label link}
\end{figure}
\pn
The number of special crossings, denoted by $\sigma$, is the number of sign changes in the Conway sequence
$ (m_1,m_2, \ldots , m_k).$
By
Proposition \ref{prop:Nlessthan},
we have $ N \le N_0 - \sigma.$
Let $ D(x)$ be the monic polynomial of degree $\sigma+ \delta$,   whose roots are the abscissas
of the $\sigma$ special crossings and the abscissas of the $\delta$ nodes
that are not crossings. The polynomial $D(x)$ is real.

Let  $ \mathcal{V}   $ be the vector space  of polynomials $ V (x,y) \in \CC [ x, y] $ such that
$$ \deg  \Bigl ( V ( x(t), y(t) ) \Bigr) \le 2b-4 .$$
The monomials $ x^{\alpha} y^{\beta}$ such that $ 3 \alpha + b \beta  \le 2b-4 $ form a basis of
$ \mathcal{V} ,  $ and it is not difficult to see that the number of these monomials is $b-1.$

Let $\mathcal{F}$ be the vector space of complex functions defined on the set of nodes of $C.$
The restriction induces a linear mapping
$\iota: \mathcal{V}  \rightarrow   \mathcal{F}   $
between spaces of the same dimension. If $U(x,y) $ is in the kernel of
$\iota$,
then we have
$ U(x(t), y(t) ) =0$ for $2b-2$ values of $t.$ Since $ \deg U( x(t),
y(t) ) \le 2b-4$,
we see that
$U(x,y)=0$.
Hence $\iota$ is an injective mapping
and then it is  an isomorphism.

For each non-special crossing with parameters $ (t_i, s_i)$, let $h_i$ be a real number
in the open interval $ (z(t_i), z(s_i))$.
Since $\iota$ is an isomorphism,
there exists a unique polynomial $ V(x,y)$ such that
$ V(x_i, y_i) = h_i D(x_i) $ for each non-special crossing $( x_i, y_i)$,
and $V(x,y)=0$ for all other nodes of $C.$
By uniqueness, we see that $V(x,y)$ is a real polynomial.
Let us consider the rational function $h(t)$ defined by
$$
h(t)= \Frac {  V(x(t), y(t)) }{ D(x(t)) }.
$$
 Each parameter $t$ of a special crossing (or  \emph{special parameter})
 is a zero of the numerator and a simple zero of the denominator.
Consequently, the function $h(t)$ is defined for all crossing parameters.
 Up to perturbing $z(t)$ by a constant if necessary, we can assume
 that $ z(t_i) \ne h(t_i) $ for all crossing parameters $t_i$.

Now, we shall prove that the  polynomial equation
\[
z(t) D(x(t)) - V(x(t),y(t)) =0
\label{eq:DV}
\] has at least  $2b-3$ distinct roots.

First, the two parameters $t,s$ of a  node  such that $ V(x,y)=D(x)=0 $ are roots of this equation.
The number of such roots is  $2 (\sigma+ \delta).$
The other roots are the zeroes of the rational function $ \Delta (t) = z(t) - h(t).$
\pn
An interval $[r,s]\subset \RR$ is called \emph{minimal} if $r,s$
are two non-special node parameters, and if $s>r$ is minimal for this property.
In other words,
there is no non-special node parameter $\tau$ in $(r,s)$.
The number of minimal intervals is
exactly $2(N_0-\sigma)-1$.

We claim that every minimal interval contains a zero of $ \Delta(t)$
that is not a node parameter.
Then the number of distinct roots of Equation (\ref{eq:DV})  must be at least
$$2(N_0 -\sigma)-1 +2(\sigma+ \delta) = 2 (N_0+\delta) -1 = 2b-3, $$
and the degree of the equation must be at least $2b-3.$

\begin{figure}[!ht]
\begin{center}
\psfrag{A}{$A$}
\psfrag{B}{$B$}
\psfrag{C}{$C$}
\psfrag{a}{}
\begin{tabular}{cccc}
\imagetop{\scalebox{.8}{\includegraphics{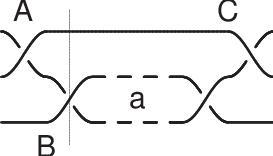}}}&\quad
\imagetop{\scalebox{.8}{\includegraphics{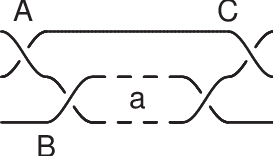}}}&\quad
\imagetop{\scalebox{.8}{\includegraphics{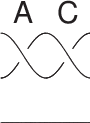}}}&\quad
\imagetop{\scalebox{.8}{\includegraphics{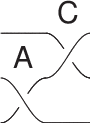}}}
\end{tabular}
\end{center}
\vspace{-.3cm}
\caption{\small The rightmost sub-arc $AC$ (ordinary cases)}
\label{fig:rightmost-o}
\end{figure}

\begin{figure}[!ht]
\begin{center}
\psfrag{A}{$A$}
\psfrag{B}{$B$}
\psfrag{C}{$C$}
\psfrag{a}{$m_k$}
\begin{tabular}{cc}
\imagetop{\scalebox{.8}{\includegraphics{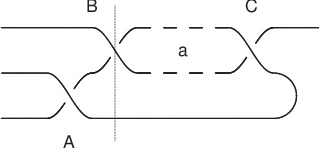}}}&\quad
\imagetop{\scalebox{.8}{\includegraphics{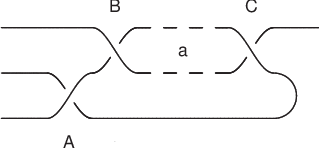}}}
\end{tabular}
\end{center}
\vspace{-.3cm}
\caption{\small The rightmost sub-arc $AC$ (exceptional case)}
\label{fig:rightmost-e}
\end{figure}

Since $ \deg V(x(t), y(t) ) \le 2b-4 $, we deduce that
$$ \deg ( z(t) D(x(t)) = c+3 (\delta +\sigma) \ge 2b-3,$$
and then
$ b+c \ge 3 (b-1-\delta -\sigma  ) = 3 (N_0-\sigma) \ge 3N,$
which conclude the proof in this case.
\pn
Let us prove our  claim. To do so, we study the sign of the rational function $\Delta (t)$
on the minimal interval $[r,s]$.
Let $j$ be the number of  special parameters contained in $[r,s]$,
and let $t_0=r$, $t_{j+1} = s$.
If $j \ne 0,$ then let $t_1 < t_2 < \ldots < t_j$ be the special
parameters contained in $  [r,s].$
The function $\Delta (t)$  is defined  for each $t_i$,
and we have $ \Delta (t_i) \ne 0 .$
The poles occur for the parameters $\tau \in [r,s]$ such
that $ D ( x(\tau) )=0$
and $(x(\tau), y (\tau) )$ is not a crossing, they are simple poles.
Let  $ [t_h, t_{h+1}] $ be the interval where the function
$ x(t), \  \tau \in [r, s] $
has a maximum. On this interval there is either one pole and no alternation
in the Gauss sequence of the knot,
or no pole and one alternation.
Figures \ref{fig:rightmost-o} and \ref{fig:rightmost-e} shows the main  cases,
the interval $ [t_h , t_{h+1} ] $ corresponds to
the rightmost sub-arc $AC$ of the arc  parametrised by $ [r,s]$.

\begin{figure}[!ht]
\begin{center}
\psfrag{A}{$A$}
\psfrag{B}{$B$}
\psfrag{C}{$C$}
\psfrag{a}{}
\begin{tabular}{ccc}
\imagetop{\scalebox{.8}{\includegraphics{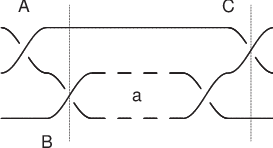}}}&\quad
\imagetop{\scalebox{.8}{\includegraphics{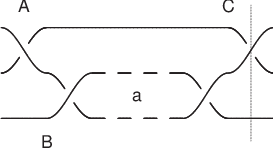}}}&\quad
\imagetop{\scalebox{.8}{\includegraphics{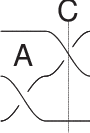}}}
\end{tabular}
%\vspace{-.3cm}
\caption{\small The other sub-arcs $AC$}
\label{fig:sub-arc}
\end{center}
\end{figure}

\pn
On the other  intervals $ [ t_{i}, t_{i+1}], \  i \ne h $  there
is either one pole and one alternation, or no pole and no alternation, see
Figure \ref{fig:sub-arc}.
Consequently, we see that $\Delta (r) \Delta (s)  <0$ if and only if the number of poles
contained in $ [r,s]$ is even. On the other hand, the number of sign changes in $[r,s]$ of the function
$\Delta(t)$ is odd if and only if $\Delta (r) \Delta (s)  <0$.
Consequently, whatever the sign of $\Delta (r) \Delta (s)$ may be, there must be at least one
$u \in [r,s]$ which is not a pole, and where $\sign{\Delta(t)}$ changes. Hence, $u$ is a root of Equation (\ref{eq:DV}), which proves the claim.
\pn
In the general case, there may be some $m_i=0$ in the diagram
$C(m_1, m_2, \ldots, m_k)$,  where $2 \leq i \leq k-1$.
We shall inductively select some relevant crossings, and ignore the others.
\pn
If there is a subsequence of the form $ ( m, 0, -n), m\ge n >0$, then we declare the last $2n$ crossings
irrelevant, and we consider the new Conway sequence where $ (m,0,-n)$ has been changed to $ m-n.$
We iterate this selection (by elimination)
until we obtain a diagram $ C(m_1, \ldots , m_k)$ such that  $m_i \ne 0$, for $i=2, \ldots k-1.$
Then,  considering  only the relevant crossings, we choose the special crossings.
We also define $D(x)=V(x,y) =0$ for the special crossings,
the irrelevant crossings and the nodes that are not crossings.
The rest of the proof is similar to the preceding one, except that the number of poles on each
minimal interval $[r,s]$ may be increased by an even number, which does not change the
sign of $\Delta(r) \Delta(s)$.
\EPf
\pn
In \cite{KP4},
it is proved that every two-bridge knot of crossing number $N$ admits an
explicit parametrisation of the form $(T_3, T_{b}, C) $
where $T_n$ is the Chebyshev polynomial of degree $n$ defined by
$T_n( \cos t) = \cos nt$, and $b+ \deg C =3N$.
Moreover, the harmonic knot
$\H(3,b,c) : (T_3,T_b,T_c)$, where $b<c<2b$, $b+c\equiv 0 \Mod 3$
has crossing number $N= \frac 13 (b+c)$, see \cite[Corollary~6.6]{KP4}.
\pn
Combining with Theorem \ref{thm:b+c}
we deduce the following.
\begin{corollary}\label{prop:lowup}
The lexicographic degree $(3,b,c)$ of a two-bridge knot of crossing number $N$ satisfies:
%Let $K$ be a two-bridge knot with crossing number $N$.
%Then the lexicographic degree $(3,b,c)$ of $K$ satisfies:
\begin{align*}
%%pvk 3<b<c, \, b \not \equiv 0 \Mod 3,\, c \not \in 3 \NN + b \NN,\, b+c \geq 3N, \\
3<b<c<2b \, b \not \equiv 0 \Mod 3,\, b+c \equiv 0 \Mod 3,\, b+c \geq 3N, \\
(3,N+1,2N-1) \leq (3,b,c) \leq (3,\pent{3N-1}2,\pent{3N}2 +1).
\end{align*}
Moreover, these inequalities are best possible.
\end{corollary}
\Pf
The transformation $(x,y,z) \mapsto (x,y-\lambda x^u, z- \mu x^v y^w)$, where
$u,v,w$ are nonnegative integers  %$\in \NN$
and $\lambda, \mu \in \RR$, does not
change the nature of the knot. This ensures that
$b \not \equiv 0 \Mod 3$ and $b+c \equiv 0 \Mod 3$.
Next, it is proved in \cite{KP4},
that every two-bridge knot admits a polynomial parametrisation of lexicographic degree
$(3,b,c)$, with  $b+c=3N$.
This implies that $b \leq \pent{3N-1}2$.
Furthermore if $b=\pent{3N-1}2$, then
$c \leq 3N-\pent{3N-1}2 = \pent{3N}2 +1$.
If $\gamma:\RR\to\RR^3$ is a polynomial parametrisation of degree
$(3,b,c)$ of a knot, then by forgetting the last coordinate we obtain a polynomial
map $\RR\to\RR^2$ of degree $(3,b)$ with at least $N$ crossings. The
genus formula implies that $b\ge N+1$. In the case $b= N+1$,
Theorem \ref{thm:b+c} implies that $c\ge 2N-1$.
\pn
Let us show that these bounds are sharp.
If $N \not \equiv -1 \Mod 3$, then the harmonic knot
$\H(3,N+1,2N-1)$ is of degree $(3,N+1,2N-1)$.
If $N \equiv -1 \Mod 3$, then $b \geq N+2$ and then $c \geq 2N-2$.
In this case, the harmonic knot $\H(3,N+2,2N-2)$ is of degree $(3,N+2,2N-2)$.
The twist knots of crossing number $N$ are of maximal degree $(3,\pent{3N-1}2,\pent{3N}2 +1)$, see \cite{BKP2}.
\EPf
\begin{remark}
%Unfortunately, the
The degree of a harmonic knot may be smaller than the degree
of its harmonic diagram. For example the knot
$\H(3,11,16) = \overline{9}_{17}$ is of degree $(3,10,17)$,
see Table \ref{table:9}.
\end{remark}

\begin{proposition}\label{prop:b+c=3N}
Let $(3,b, c)$ be the lexicographic degree of
a two-bridge knot of crossing number $N$.
If $b \leq N+3$ or $b = \pent{3N-1}2$ then we have  $b+c = 3N$.
\end{proposition}
\Pf
By Theorem \ref{thm:b+c}, we have $b+c \geq 3N$, and $b+c = 3N$ if
$b= \pent{3N-1}2$ by Corollary  \ref{prop:lowup}.
Hence we assume now that $b \le N+3$.
Let $\gamma(t) = (x(t),y(t),z(t))$ be a polynomial representation of
our knot $K$ of degree $(3,b,c)$, and denote by
$D=C(m_1, \ldots, m_k)$ the trigonal diagram of $\gamma$.
If $s$ denote the number of sign changes in the Gauss sequence of the
parametrisation $\gamma$, we clearly have $c\leq s$. Hence it remains us to obtain an
upper bound for $s$, using Propositions \ref{prop:Nlessthan} and \ref{prop:s}.

Let $N_0=\sum \abs{m_i}$, and
$\sigma$ be the number of sign changes in the sequence
$(m_1, \ldots, m_k)$.
Combining  Propositions \ref{prop:Nlessthan} and the genus formula for
plane curves, we obtain
\begin{equation}\label{eq1}
N+\sigma  \leq N_0 \leq b-1.
\end{equation}
First, suppose that $b=N+3$.
In this case $N \not \equiv 0 \Mod 3$, and
$c \equiv -N \Mod 3$, by Corollary \ref{prop:lowup}. Consequently $c \not \equiv 2N-1 \Mod 3$ and $c \not \equiv 2N-2 \Mod 3$.
%pvk both $2N-1$ and $2N-2$ belong to $\langle 3,N+3 \rangle$.
Hence we only have to prove that $c \leq 2N-1$.
\begin{enumerate}
\item First, suppose  that $D=C(x,m,0,-n,-y)$
  with $mn>0$, see Figure \ref{fig:m0n-1}. Since $N_0\le N+2$, we
  necessarily have $\abs m=1$ or $\abs n=1$.   Without loss of generality,  we can
  assume that $n=1$ and $m>0$.
  Consider the diagram
  $D'=C(x,m-1,-y)$
  obtained by a type-II Reidemeister move on $D$.
  The diagram $D'$ has
$N_0-2=N$
  crossings, and then is an alternating
diagram of $K$. Consequently  the number $s'$ of sign changes
in the Gauss sequence of $D'$ is $s'=2N-1$.

If $(x,m) \not = (1)$ and $(n,y) \not = (1)$, then we have $s=s'=2N-1$ and consequently
$c \leq 2N-1$, see Figure \ref{fig:m0n-1}.
\begin{figure}[!ht]
\psfrag{y}{}
\psfrag{b}{\!\!\!\!$m-1$}
\begin{center}
\begin{tabular}{ccc}
{\scalebox{.5}{\includegraphics{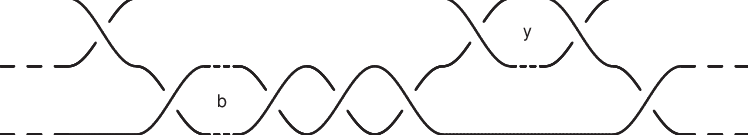}}}&\quad&
{\scalebox{.5}{\includegraphics{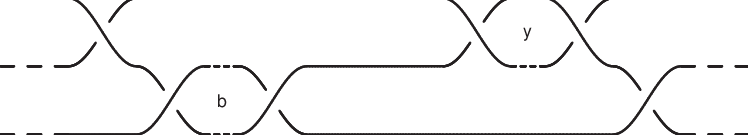}}}
\end{tabular}
\end{center}
\vspace{-0.3cm}
\caption{\small  $C(x,m,0,-1,-y)  \mapsto C(x,m-1,-y)$}
\label{fig:m0n-1}
\end{figure}
\pn
If we have $(x,m)=(1)$ or $(n,y)=1$, then we can suppose $(n,y)=(1)$
and $D=C(x,m,0,-1)$.
If we change the nature of the last two crossings,
then we obtain another diagram
$\tilde D = C(u,m-1,0,-1,0,1)$ of
$K$ with the same $xy$-projection.
By the
previous case, we see that the number of sign changes in the Gauss
sequence of $\tilde D$ is
$\tilde s = 2N-1$. Consequently $\tilde D$ is of degree at most $(3,N+3,2N-3)$.
\pn
\item Then,
suppose that $D=C(x,n,-1)$
(the case $D=C(1,-m,y)$ is similar).
By changing the nature of the last two crossings of $D$, we obtain
another diagram $\tilde D=C(x,n-1,0,-1,1)$ of the same knot,
see Figure \ref{fig:1-m}.
By
case $1$ above,
we see that the number of sign changes in the Gauss sequence of
$\tilde D$ is $\tilde s=2N-1$ and we deduce that $\tilde D$ is of degree
$(3,N+3,2N-3)$.
\begin{figure}[!ht]
\psfrag{a}{}
\psfrag{b}{\!\!\!\!$m-1$}
\begin{center}
\begin{tabular}{ccc}
{\scalebox{.5}{\includegraphics{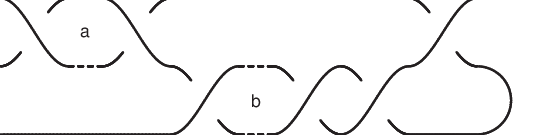}}}&\quad&
{\scalebox{.5}{\includegraphics{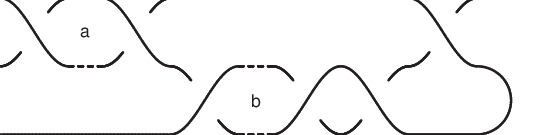}}}
\end{tabular}
\end{center}
\vspace{-0.3cm}
\caption{\small  $C(x,m,-1) = C(x,m-1,0,1,-1) \mapsto \tilde D = C(x,m-1,0,-1,1)$}
\label{fig:1-m}
\end{figure}
\item  Now, suppose that $D$ is not in the cases $1$ and $2$ above.
   If $\sigma=2$, then $N_0=N+2$ and $\sigma_2 \leq 1$. By Proposition
   \ref{prop:s}, we obtain
   $s=(2N_0-1) - 3 \sigma + 2 \sigma_2 \leq 2N-1$.

   If $\sigma<2$, then
   by Lemma \ref{rem:N<N0} we have $m_1 \cdot m_k=0$.
 Consider the reduced diagram $D'$. If $\sigma=0$, then
    $D'$ is alternating and has $N_0=N$ crossings. Its Gauss sequence is alternating and has $s'= 2N-1$ sign changes.
    If $\sigma=1$, then $D'$ may have $N'_0=N$ or $N'_0=N+1$
    crossings. If $N'_0=N$ then $D'$ is alternating and there are
    $s'=2N-1$ sign changes in its Gauss sequence. If $N'_0=N+1$, then
    $D'$ is not alternating and $\sigma'_1=1$. We thus have $s'= 2N+1
    - 3 =2N-2$
    by  Proposition \ref{prop:s}.

    We then choose a polynomial of degree $c \leq s' \leq 2N-2$ as an height function for the reduced diagram $D'$. If $m_1=0$ (resp. $m_k=0$), the signs of the $\abs{m_2}$ (resp. $\abs{m_{k-1}}$)  crossings do not affect the nature of the knot.
\end{enumerate}
At the end we find a polynomial height function $z(t)$ of degree $c \leq 2N-1$.
\pn
If $b=N+2$, then $N \not \equiv 1 \Mod 3$ and $2N-1 \in \langle 3,b
\rangle$. Hence again,
we only have  to prove  $c \leq 2N-1$.
By Inequality $(\ref{eq1})$, we may have $N_0=N$ or $N_0=N+1$.
\begin{enumerate}
\item
If $N_0=N$, then the diagram is alternating and $s \leq 2N-1$.
\item
  If $N_0=N+1$, then $\sigma \leq 1$.
  If $\sigma=1$, then $s\le 2N-1$ by Proposition \ref{prop:s}.
If $\sigma=0$, then $m_1 \cdot m_k =0$
by Lemma \ref{rem:N<N0}.
The reduced diagram is alternating and its Gauss sequence has $s'\leq
2N-1$ sign changes
and so $c \leq 2N-1$.
\end{enumerate}
At the end we find a polynomial
function $z(t)$ of degree $c \leq 2N-1$.
\pn
If $b=N+1$ then $N_0=N$ and the diagram is alternating. We thus have
$c \leq s \leq 2 N -1$.
\EPf
\pn
We deduce
\begin{proposition}\label{prop:N<=11}
The lexicographic degree $(3,b,c)$ of a knot with crossing number
$N\le 11$ satisfies $b+c=3N$.
\end{proposition}
\Pf
By Proposition \ref{prop:lowup}, we have
$(3,b,c) \leq (3,\pent{3N-1}{2},\pent{3N}2+1)$.
If $b\leq N+3$ or $b = \pent{3N-1}{2}$, we conclude using Proposition \ref{prop:b+c=3N}.
If  $b\ge N+4$ and $b < \pent{3N-1}{2}$, then $N=11$, and $b=15$ which is impossible since $b$ is not divisible by 3.
\EPf

\section{Degrees of trigonal plane diagrams}\label{sec:plane}

Thanks to the relation $b+c=3N$ established in Proposition \ref{prop:N<=11},
we are now reduced to study plane trigonal curves.
It is enough to determine the smallest integer $b$ such that the
$xy$-projection of some diagram of $K$
admits a polynomial parametrisation of degree $(3,b)$.
%pvk to determine the lexicographic degree of a two-bridge knot $K$ with crossing number at most 11.

Given a long knot diagram $D$ in $\RR^3$, we denote by $\abs D$
its projection to $\RR^2$ (i.e. we forget about the sign of the crossings).
If $D=C(m_1,\ldots, m_k)$, we use the
notation $|D|=D(|m_1|,\ldots, |m_k|)$.
An isotopy of $\RR^2$ is called an {\em $\mathcal L$-isotopy} if it
commutes with the projection
$\RR^2\to\RR$ forgetting the second coordinate.
\begin{definition}
The \emph{algebraic degree} of
$\abs D$ is the minimal integer $b$ such that there exists a
real algebraic curve
$\gamma:\CC\to\CC^2$ of bidegree $(3,b)$ such that $\gamma(\RR)$ is
$\mathcal L$-isotopic to $\abs D$.
\end{definition}
We first establish a lower bound for polynomial curves in Proposition
\ref{prop:Bezout}.

\subsection{Lower bounds on degrees of plane trigonal diagrams}

\begin{proposition} \label{prop:Bezout}
Let $|D|$ be the plane diagram
$D(m_1, m_2, \ldots , m_k), $  with $ \  m_i \ge 2$ for  $i= 1,\ldots, k .$
Then the
algebraic
degree of $|D|$ is at least
$3k-1$.
If in addition we have $ m_i \ge 3 $ for some $i$, then the
algebraic
degree of $ |D|$ is at least
$3k+1$.
\end{proposition}
\Pf
Let $ \gamma(t): (x(t), y(t)) $ be a polynomial parametrisation of
$|D|$ with $x(t)$ of degree 3, and let $C$ be the image of $\gamma$.
The complement of $C$
contains
$ m_j -1 $
disks
corresponding to the $j$th group of crossings of $|D|$.
Let us choose a point $P_j$ in one of these
disks.
There is a polynomial curve of equation $ y= P (x)$ with $\deg P= k-1$
containing the $k$ points $P_j$.
\begin{figure}[!ht]
\begin{center}
 {\scalebox{.6}{\includegraphics{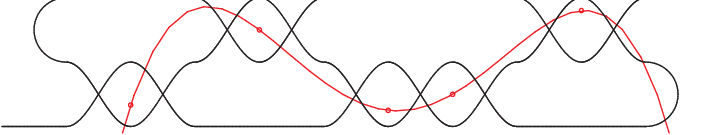}}}
\end{center}
\vspace{-10pt}
\caption{\small  The plane diagram $D(2,2,3,2)$}
\label{fig:9_15_alt}
\end{figure}
Since the number of intersections of this curve and $C$ is at least
 $ 2k + (k-1) = 3k-1 > 3(k-1)$, we deduce that $\deg(y(t))  \ge  3k-1$.

If in addition some $m_i \ge 3,$ we choose one more
point $P_{k+1}$ in another
disk
of the $i$th group of two-sided domains.
Then we count the intersections of $C$ with a curve $ y = P(x) $
$\deg (P(x))=k$ containing the $k+1$ points $P_j, \ j= 1,\ldots,k+1$.
Since this number is at least $ 2 ( k+1) + (k-1) = 3k+1 > 3k, $ we deduce that
$\deg (y(t)) \ge 3k+1$ (see Figure \ref{fig:9_15_alt} in the
case of $D(2,2,3,2)$).
\EPf
\subsection{Application of Orevkov's braid theoretical method}
To obtain lower bounds on the algebraic degree $b$, it is convenient to enlarge
the category of objects under interest, and to consider
\emph{real pseudoholomorphic curves} rather than real algebraic
curves. Doing so, we can use the
full power of
 the braid theoretical approach developed by Orevkov to
study real  curves in
$\CC^2$. Using this strategy,
we determined
in \cite{BKP2} the lexicographic degree of all torus knots
$C(m)$ and generalised twist-knots $C(m,n)$.
We refer to \cite[Section 3.2]{BKP2} for the definition of  a
real pseudoholomorphic curve $\gamma:\CC\to\CC^2$ of
bidegree $(3,b)$ where $b$ is a positive integer.
%pvk $\in\ZZ_{\ge 1}$.
Recall that a real
algebraic map $\gamma:  \CC \to \CC^2$ of degree $(3,b)$
is  an example of a real
pseudoholomorphic curve of bidegree $(3,b)$.
Without loss of generality, we only consider in this text nodal
pseudoholomorphic curves.

\begin{definition}
The \emph{pseudoholomorphic} degree of $\abs D$ is the
minimal integer $b$ such that there exists a real pseudoholomorphic
curve $\gamma:\CC\to\CC^2$ of bidegree $(3,b)$ such that $\gamma(\RR)$ is
$\mathcal L$-isotopic to $\abs D$. It is not greater than the algebraic degree of $\abs D$.
\end{definition}
Recall  that the {\em group of braids with $3$-strings} is defined as
$$B_3=\langle \sigma_1, \sigma_{2}|
\sigma_1\sigma_{2}\sigma_1=\sigma_{2}\sigma_{1}\sigma_{2}\rangle.$$
We  refer to  \cite[Sections 2 and 3]{BKP2} for an algorithm that
associates
an \emph{$\mathcal L$-scheme} and
a braid $\bb_C\in B_3$ to any real pseudoholomorphic curve
$C=\gamma(\CC)$, with  $\gamma:\CC\to\CC^2$   a
real pseudoholomorphic curve
of
bidegree $(3,b)$. A braid $\bb \in B_3$ is said to be
\emph{quasipositive} if it can be written in the form
\[
\bb=\prod_{i=1}^l w_i \sigma_1 w_i^{-1} \qquad \mbox{with
}w_1,\cdots,w_l\in B_3.
\]
Note that
a braid with algebraic length 0 is quasipositive if and only if
it is the trivial braid. The
quasipositivity problem in $B_3$ has been solved by Orevkov \cite{O3}.
We will use the following proposition in order to obtain lower bounds
in lexicographic degree of knots.
\begin{proposition}\label{prop:braid}
 Let   $\gamma:\CC\to\CC^2$  be a
real pseudoholomorphic curve of
bidegree $(3,b)$, and let $C=\gamma(\CC)$. We denote by $\pi:\CC^2\to
\CC$ the projection to the first coordinate, and we assume that
the two critical points of the map
$\pi\circ \gamma$ are real. Then the braid $\bb_C$
satisfies the three following properties:
\bn
\item[(i)] $\bb_C$ is quasipositive;
\item[(ii)] the closure of $\bb_C$ is a link with three components;
\item[(iii)] the linking number of any two strings of $\bb_C$ is non-negative.
\en
\end{proposition}
\Pf
  Property (i) is a consequence of \cite[Proposition 7.1]{O1}. Properties (ii) and
  (iii) are easy consequences of the Riemann--Hurwitz formula applied to
  the map $\pi\circ \gamma$, see \cite[second paragraph of the proof
    of Proposition 3.1]{BKP2}.
\EPf
\begin{remark}\label{rem:square quasi}
Proposition \ref{prop:braid}(i) can be strengthened in order to get an
equivalence. Given $l\in\ZZ_{\ge 1}$ and $\varepsilon=1$ or $2$,
we define $b=3l-\varepsilon$.
Let $\LL_S$ be the trigonal $\LL$-scheme
\[
 \mathcal L_S= \bullet_{i_1} \cdots \bullet_{i_{\alpha_1}} \subset_{*} \times_{j_1}
\cdots \times_{j_{N_0}} \supset_{*} \bullet_{i_{\alpha_1+1}} \cdots \bullet_{i_{\alpha}} \star,
\]
where $\star=\downarrow$ or $\uparrow$ if $\varepsilon=1$, and
$\star=\vee$ or $\land$ if $\varepsilon=2$.
Then one can associate a braid $\bb_C$, depending on $b$,
to the $\mathcal L$-scheme
$\mathcal L_S$ using the algorithm given in \cite[Section
  2.2]{BKP2}. Following  \cite[Proposition 7.1]{O1}, we have that
$\mathcal L_S$  is realised by a real pseudoholomorphic curve of
bidegree $(3,b)$ in $\CC^2$ if and only if the braid $\bb_C$ can be
written in the form
\[
\bb_C=\prod_{i=1}^{\ell} w_i \sigma_1^2 w_i^{-1} \qquad \mbox{with
}w_1,\ldots,w_{\ell}\in B_3.
\]
Note that in this case, we necessarily have
$\deg \bb_C = 2 \ell =b-1-\alpha - N_0$.
\end{remark}
\begin{remark}
Proposition \ref{prop:Bezout} also holds for the pseudoholomorphic
degree of a plane trigonal diagram, and the proof is essentially the
same. Nevertheless we will not need this more general version here.
\end{remark}
We end this section by proving a slight generalisation of
\cite[Proposition 3.1]{BKP2}.
\begin{proposition}\label{prop:even}
Let $D=C(m_1, \ldots, m_{k})$ be a trigonal diagram of a knot $K$, with
$m_1, \ldots,$ $m_{k-1}$ even integers.
As usual, we
define $N_0= m_1 + \cdots + m_k$.
If $\gamma:\CC\to\CC^2$ is a  real rational pseudoholomorphic curve
of bidegree $(3,b)$ such
that $\gamma(\RR)$ is $\LL$-isotopic to
$|D|$,
then
$2 b \geq 3N_0 -2$.
\end{proposition}
\Pf
Let us write $b=3l-1$ or $b=3l-2$,
let $\alpha$ be the number of solitary nodes of $C=\gamma(\CC)$,
and
$2\beta$ be the number of complex conjugated nodes.
By the genus formula, we have $$N_0+\alpha+2\beta = b-1.$$

The $\LL$-scheme realised by $C$ has the form
$$
 \bullet_{i_1} \cdots \bullet_{i_{\alpha_1}} \subset_{*} (\times_{j_1})^{m_1}
\cdots (\times_{j_k})^{m_k} \supset_{*} \bullet_{i_{\alpha_1+1}} \cdots \bullet_{i_{\alpha}
} \star,
$$
where $\star=\downarrow,\uparrow, \vee$ or $\land$. The braid $\bb_C$
has 3 components $L_1$, $L_2$ and $L_3$, and
$\lk(L_i,L_j)\geq 0$ by Proposition \ref{prop:braid}.
Furthermore, as in \cite[proof of Proposition 3.1]{BKP2}, we have
$0 \leq \lk(L_i,L_j)\leq \beta$.

By the assumptions made on $D$,
there are two strings of $\bb_C$, say $L_1$ and $L_3$,
that do not cross at the crossing points of $\RR C$.
Each $\bullet_{j} \bullet_{j'}$ contributes at least $-1$ to
$\lk(L_1,L_3)$.
Hence as in \cite[Proof of Proposition 3.1]{BKP2},
we obtain
$$2 \beta \geq 2 \lk(L_1,L_3) \geq l - \alpha -2,$$
and thus
$$ b-1 = N_0 + \alpha + 2 \beta \geq N_0+l-2.$$
We then deduce
$3b -3N_0 \geq 3l-3 \geq b-2$, and  $2b \geq 3N_0 -2$.
\EPf

\subsection{The $\TT$-reduction}\label{sec:upper}
\begin{definition}\label{def:r4}
Let $x,y$ be (possibly empty) sequences of nonnegative integers
and $m,n$ be nonnegative integers.
The plane diagram
$D( x, m,n,y)$ is called a
\emph{$\TT$-reduction} of the  diagram
$D(x, m+1,1,n+1, y) $ (see Figure \ref{fig:T}).
\begin{figure}[!ht]
\begin{center}
\psfrag{a}{$m-1$}
\psfrag{b}{$n-1$}
{\scalebox{.6}{\includegraphics{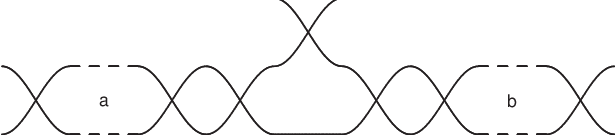}}} \quad\quad
\psfrag{a}{$m$}
\psfrag{b}{$n$}
{\scalebox{.6}{\includegraphics{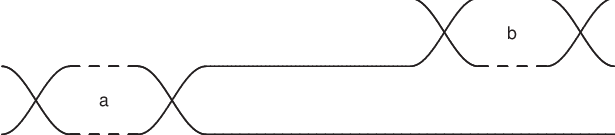}}}
\end{center}
\caption{\small $\TT$-reduction}
\label{fig:T}
\end{figure}
\end{definition}
\pn
Propositions \ref{prop:flip} and \ref{prop:perturb} below relate the
pseudoholomorphic and algebraic degrees of  two plane trigonal diagrams
differing by a $\TT$-reduction.
\begin{proposition}\label{prop:flip}
Let $\abs {D_1}$ and $\abs {D_2}$ be two plane trigonal diagrams such that
$\abs {D_2}$ is obtained from $\abs {D_1}$ by a $\TT$-reduction.
If $\abs {D_1}$ has pseudoholomorphic degree $b$,
then $\abs {D_2}$   has  pseudoholomorphic degree $b-3$.
\end{proposition}
\Pf
Let $\abs{D_1}=D(m_1,\ldots,m_k)$ and
$\abs{D_2}=D(n_1,\ldots,n_l)$.
Suppose that there exists a real pseudoholomorphic
curve $\gamma_1:\CC\to\CC^2$ of bidegree $(3,b)$ such that $\gamma_1(\RR)$
is $\mathcal L$-isotopic to $\abs{D_1}$, and suppose that its
associated $\mathcal L$-scheme is
\[
 \bullet_{i_1} \cdots \bullet_{i_{\alpha_1}} \subset_{*} (\times_{j_1})^{m_1}
\cdots (\times_{j_k})^{m_k} \supset_{*} \bullet_{i_{\alpha_1+1}} \cdots \bullet_{i_{\alpha}} \star.
\]
The braid associated to $\gamma_1$ is the same that the braid
associated to the $\mathcal L$-scheme
 \[
 \bullet_{i_1} \cdots \bullet_{i_{\alpha_1}} \subset_{*} (\times_{j_1})^{n_1}
\cdots (\times_{j_l})^{n_l} \supset_{*} \bullet_{i_{\alpha_1+1}} \cdots \bullet_{i_{\alpha}} \star.
\]
Hence according to Remark \ref{rem:square quasi},
there exists a real pseudoholomorphic curve
$\gamma_2:\CC\to\CC^2$ of bidegree $(3,b-3)$ such that $\gamma_2(\RR)$ is
 $\mathcal L$-isotopic to $\abs{D_2}$.
\EPf
\begin{corollary}\label{cor:lexdeg1}
The pseudoholomorphic degree of the plane diagram $D(0,n)$ is
$\pent{3n}2+1$.
\end{corollary}
\Pf
The plane diagram $D(0,n)$ is obtained by a $\TT$-reduction
from
$D(1,1,n+1)$. Since
$D(1,1,n+1)$ and $D(2,n)$
may be reduced to each other by slide isotopies, they have the same
pseudoholomorphic degree
by Proposition \ref{cor:isotopy2}.
By \cite[Theorem 3.9]{BKP2}, the degrees are
$\pent{3n}2+4$, which completes the proof.
\EPf
\subsection{The $\TT$-augmentation}\label{sec:augmentation}
Proposition \ref{prop:flip} admits a weaker version for
the algebraic degree of a plane diagram. We make use
the $\TT$-augmentation that consists in adding a triangle of crossing points in a given plane diagram.
\begin{proposition}\label{prop:perturb}
Let $\abs{D_1}$ and $\abs{D_2}$ be two plane trigonal diagrams
such that $\abs{D_2}$ is
obtained from $\abs{D_1}$ by a reduction $\TT$.
If $\abs{D_2}$ has algebraic degree $b-3$, then $\abs{D_1}$
has algebraic degree at most $b$. Furthermore, if the pseudoholomorphic
degree of $|D_2|$ is also $b-3$, then $\abs{D_1}$
has algebraic degree exactly $b$.
\end{proposition}
%\Rpvk{Ne faut-il pas etre plus explicite : "If $\abs{D_2}$ admits a polynomial parametrisation of degree $(3,b-3)$ then $\abs{D_1}$ admits a polynomial parametrisation of degree $(3,b)$".}
\Pf
The last assertion simply follows from the fact that a real rational
algebraic curve in $\CC^2$ is a pseudoholomorphic curve.
Let
\[
\begin{array}{cccc}
\gamma: & \CC &\longrightarrow & \CC^2
\\ & t &\longmapsto & (P(t), Q(t))
\end{array}
\]
be a real algebraic map with $P(t)$ of degree $3$ and $Q(t)$ of degree
$b-3$, and such that $\gamma(\RR)$
 is
$\LL$-isotopic to the plane diagram
$D(x,m,n,y)$, where $x,y$ are (possibly empty) sequences of
nonnegative integers and $m,n$ are nonnegative integers.
Without loss of generality, we can suppose that the line $x=0$ separates the
$m$ crossings from the $n$ crossings.
The curve parametrised by $t\mapsto (P(t),P(t)\cdot Q(t))$
has the same double points as $\gamma(\RR)$ and an additional
ordinary triple point at $(0,0)$.
For
$\eps$ small enough the curve
$(P(t+\eps), P(t)\cdot Q(t))$ is $\LL$-isotopic to either
$D(u,m+1,1,n+1,v)$ or $D(u,m,1,1,1,n,v)$, depending on the sign of
$\eps$
(see Figure \ref{exa:d3}).
\begin{figure}[!ht]
\begin{center}
{\scalebox{.8}{\includegraphics{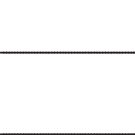}}} \quad
{\scalebox{.8}{\includegraphics{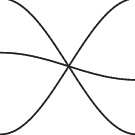}}} \quad
{\scalebox{.8}{\includegraphics{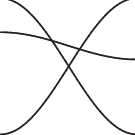}}} \quad
{\scalebox{.8}{\includegraphics{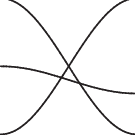}}}
\end{center}
\vspace{-10pt}
\caption{\small  Perturbation of
  a triple point in $\RR^2$}
\end{figure}
\EPf
\begin{example}\label{exa:d3}
Let us consider the polynomial parametrisation $(T_3(t),T_4(t))$ of the
diagram $D(1,1,1)$, where $T_n$ denotes the Chebyshev polynomial of
degree $n$.
We choose to add a triple point in $(-3/4,0)$,
by considering the curve $t\mapsto (T_3(t),  Q(t))$, where
$Q(t)=(T_3(t)+3/4) \cdot (T_4(t)+1)$.
\begin{figure}[!ht]
\begin{center}
    \begin{tabular}{cccc}
    {\scalebox{.9}{\includegraphics{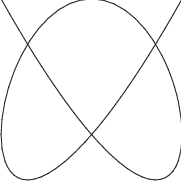}}} &
    {\scalebox{.9}{\includegraphics{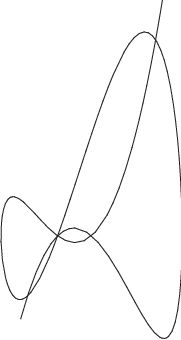}}} &
    {\scalebox{.9}{\includegraphics{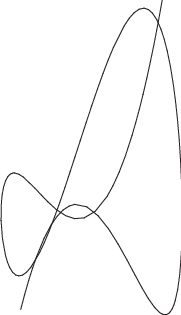}}}&
    {\scalebox{.9}{\includegraphics{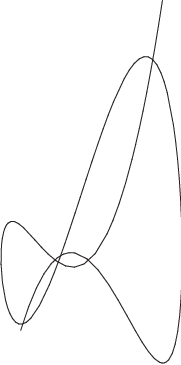}}}\\
    $(T_3,T_4)$ & $(T_3,Q)$ &  $(T_3,Q(t+\eps))$& $(T_3,Q(t-\eps))$\\
    $D(1,1,1)$ & & $D(2,1,2,1)$ & $D(1,1,1,1,1,1)$
    \end{tabular}
    \end{center}
\vspace{-10pt}
\caption{\small  Adding three crossings to the trefoil}
\label{fig:triple}
\end{figure}
Then the curve in $(P_3(t),Q(t+\eps))$ is $\LL$-isotopic
to $D(2,1,2,1)$
for $\eps >0$ small enough and is $\LL$-isotopic to
$D(1,1,1,1,1,1)$ for $\eps<0$, see Figure \ref{fig:triple}.
\end{example}

\begin{figure}[!ht]
\begin{center}
\begin{tabular}{ccccccc}
\multicolumn{3}{c}{\includegraphics[height=1.5cm, angle=0]{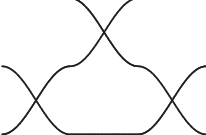}}&
\hspace{1cm}&
\multicolumn{3}{c}{\includegraphics[height=1.5cm,
    angle=0]{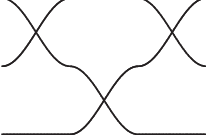}}\\
&a) & &&& b)
\\[10pt]
\includegraphics[height=1.5cm, angle=0]{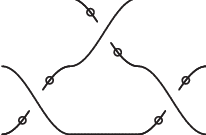}&\hspace{1cm}&
\includegraphics[height=1.5cm, angle=0]{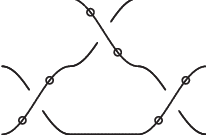}&\hspace{1cm}&
\includegraphics[height=1.5cm, angle=0]{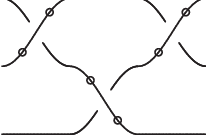}&\hspace{1cm}&
\includegraphics[height=1.5cm, angle=0]{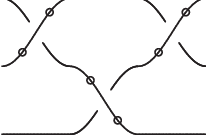}
\\ d)& &e)  &&f) && g)
\end{tabular}
\end{center}
\vspace{-10pt}
\caption{\small  Perturbation of a triple point in $\RR^3$}
\label{fig:construction}
\end{figure}

\begin{example}
  Figure \ref{fig:3T} shows that the algebraic degree of
  $D(2,2,2,1,3)$ is at most 11, starting from a
  parametrisation of the plane diagram $D(1,0)$ of degree $(3,2)$.
\begin{figure}[!ht]
\label{fig:3r3}
\begin{tabular}{cccc}
{\scalebox{.5}{\includegraphics[align=c]{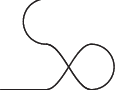}}} &
{\scalebox{.5}{\includegraphics[align=c]{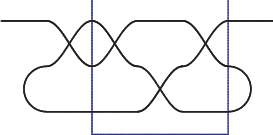}}} &
{\scalebox{.5}{\includegraphics[align=c]{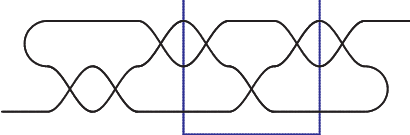}}} &
{\scalebox{.5}{\includegraphics[align=c]{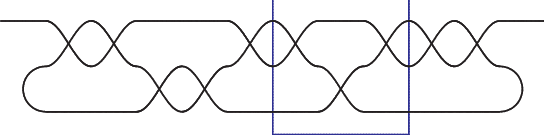}}}\\
degree $(3,2)$ &
degree $(3,5)$ &
degree $(3,8)$ &
degree $(3,11)$
\end{tabular}
\vspace{-3pt}
\caption{\small  From $D(1,0)$ to $D(2,2,2,1,3)$}
\label{fig:3T}
\end{figure}
\end{example}
Proposition \ref{prop:perturb} can be extended
to spatial trigonal curves. The next result provides
constructions of polynomial knot diagrams.
\begin{proposition}\label{prop:d+(3,6)}
Let $e=\pm 1$.
If the diagram
$C(u,m,n,v)$ has lexicographic degree $(3,b-3,c-6)$, then the diagram
$C(u,m+e,e,e+n,v)$ has lexicographic degree at most $(3,b,c)$.
\end{proposition}
\Pf
Let $t\mapsto (P(t),Q(t),R(t))$ be a parametrisation of degree
$(3,b,c)$
of the diagram  $C(u,m,n,v)$.
Up to a change of coordinates, we may assume that the part $(u,m)$
(resp. $(n,v)$) of
the diagram is contained in the half-space $x<0$ (resp. $x>0$), and
that
the three points of the diagram in the plane $x=0$ have
$z$-coordinates of the same sign.
We consider the map $\phi(t)=(P(t),P(t)Q(t),P^2(t)R(t))$. The image of
$\phi$ is a singular diagram with the three branches tangent to the
plane $z=0$ at the point $(0,0,0)$.
Extending the notations of diagram in the obvious way to this
particular case, we see that the image of
$\phi$ realises the singular diagram $C(u,m,*,n,v)$, where $*$ stands
for the triple point.
By slightly perturbing the roots of
the factor $P(t)$ of the polynomial $P(t)Q(t)$, we
obtain a polynomial $Q_1(t)$ of degree $b+3$ such that the
triple point of the curve $(P(t),P(t)Q(t))$ will be perturbed as
depicted in Figure \ref{fig:construction}a or b, depending on the perturbation
 $Q_1(t)$.
Perturbing the roots of the factor $P^2(t)$ of the polynomial
$P^2(t)R(t)$
as depicted by the blue dots on Figure  \ref{fig:construction},
we obtain a parametrisation
of the diagram whose existence is claimed in the theorem.
\EPf
\begin{example}\label{exa:d3s}
The trefoil admits the parametrisation $(T_3, T_4, T_5)$ of degree $(3,4,5)$.
We thus deduce that $6_2=C(2,1,3)$ and $6_3=C(2,1,1,2)$ admit parametrisations of degree
$(3,7,11)$. By Corollary \ref{prop:lowup},
these are the lexicographic degrees of $6_2$ and $6_3$.
\end{example}
Thanks to Proposition \ref{prop:N<=11}, we will not need
Proposition \ref{prop:d+(3,6)} to determine the lexicographic degrees of the first knots, but it 
may be useful for further results.

\section{Two-bridge knots with 11 crossings or fewer}\label{sec:lex}
Simple diagrams of two-bridge knot have been introduced in
\cite{BKP1}.
The \emph{complexity} $c(D)$ of a trigonal diagram
$D=C(m_1, \ldots, m_k)$
is defined as
$$
c(D)=k+ \Sum_{i=1}^k \abs{m_i}.
$$
\begin{definition}
We shall say that an isotopy of trigonal diagrams is a {\em slide} isotopy
if the number of crossings never increases during the isotopy, and if
all the intermediate diagrams remain trigonal.
A trigonal diagram is called a {\em simple} diagram if it
cannot be simplified into a diagram of lower complexity by using slide
isotopies only.
\end{definition}
The next two propositions motivate the consideration of simple diagrams.
\begin{proposition}[{\cite[Corollary 3.9]{BKP1}}] \label{cor:isotopy}
Let $D$ be a trigonal  diagram of a two-bridge knot.
Then by slide isotopies, it is possible to transform $D$
into a simple diagram $ C(m_1, \ldots , m_k)$
such that for $i=2, \ldots ,k$, either  $\abs{ m_i} \ne 1, \ $  or
$ m_{i-1} \, m_i  >0$.
\end{proposition}
\begin{proposition}[{\cite[Corollary 3.7]{BKP2}}]\label{cor:isotopy2}
Let ${D_1}$ and ${D_2}$ be two trigonal long knot diagrams such that ${D_2}$ is
obtained from ${D_1}$ by a slide isotopy.
Then the pseudoholomorphic degree of $|D_1|$ is greater than or equal to
the pseudoholomorphic degree of $|D_2|$.
\end{proposition}
\pn
In \cite{BKP2} we proved that the lexicographic degree of the torus knot
$C(n)$ or the twist knot $C(n,m)$ is precisely $(3,\pent{3N-1}2,\pent{3N}2+1)$ by showing first that the only simple diagrams of these knots are the
alternating diagrams and showing that the algebraic degrees of the
corresponding plane diagrams are $\pent{3N-1}2$.

\subsection{The general strategy}\label{sec:general}
Given a two-bridge knot with
crossing number $N\le 11$, our strategy to determine its lexicographic
degree consists in:
\bn
\item Find a first upper bound $b_0$ on $b$ using  constructions
from \cite{KP4} based on Chebyshev plane diagrams parametrised by
$(T_3, T_b)$, where $T_n$ is the Chebyshev polynomial $T_n(\cos t) = \cos nt$.
\item Compute all finitely many simple diagrams of $K$ with $b_0-1$
crossings or fewer. This is done by computing all continued
fractions corresponding to the Schubert fractions of $K$. % and $\overline K$.
\item For all these simple diagrams,
\bn
\item Compute a lower bound of their algebraic degree using
  Propositions \ref{prop:Bezout} and \ref{prop:even}.
\item Using $\TT$-reductions, try to  obtain  explicit constructions of
  these diagrams out of known constructions for diagrams with a lower
  number of crossings. This provides a lower bound on the
  lexicographic   degree of the knot.
\item If necessary, compute all possible braids associated
  to hypothetical plane curves of degree $b<b_0$ that are
  $\LL$-isotopic to the diagram, and check if these braids satisfy
  Proposition \ref{prop:braid}. This may improve the lower bound
  obtained in step $(a)$ above.
\item If the lower bound and the upper bound coincide, then we have
    determined the lexicographic degree of the knot.
\en
\item If the lower bound and the upper bound do not coincide,
improve the upper bound by looking at
non-simple diagrams on which one can perform
$\TT$-reductions to reduce to knots with lower crossing number.
\en
In Table \ref{table:all}, p. \pageref{table:all}, we give the lexicographic
degree of all two-bridge knots with 11 crossings or fewer.
In Tables \ref{tab6}, \ref{tab7}, \ref{tab8}, and \ref{table:9} below, we
give refinements of Table \ref{table:all} for
two-bridge knots with crossing number at most 9. The columns 1, 2 and 3 identify the knot.
The column 4 gives the lexicographic degree. The fifth column gives
the upper bound on $b$ obtained by considering Chebyshev diagrams; the
sixth column gives a diagram that can be realised in the corresponding
lexicographic degree; the last column gives the construction of the
corresponding plane diagram, when one needs to improve the bound given
by Chebyshev knots.

\subsection{Some initial diagrams}

Here we compute the algebraic degrees of a few trigonal plane
diagrams. These computation will be used in the next sections to
determine the algebraic degree of trigonal plane diagrams that
reduce to the diagrams considered in this section by $\TT$-reduction.
The next proposition is proved in \cite{KP2}.
\begin{proposition}\label{prop:torus}
The plane diagram $D(4n-1)$
has algebraic degree  $6n-2$.
\end{proposition}
This gives an explicit parametrisation for the plane diagrams
$D(3)$ and $D(7)$.
\begin{lemma}\label{lem:constr}
We give below the algebraic degree of
a few plane diagrams (see Figure \ref{fig:ddd} for the image of a
polynomial parametrisation of the given degree).
\bi
\item $b=1$: $D(0,0)$
\item $b=2$: $D(0,1)$
\item $b=4$: $D(0,2)$, $D(2,1)$
\item $b=5$: $D(0,1,1,0)$, $D(2,2)$, $D(1,1,1,1)$, $D(0,3)$, $D(1,2,0)$
\item $b=7$: $D(5)$, $D(1,4)$, $D(0,4)$
\ei
%\begin{center}
%\begin{tabular}{|c|l|}
%\hline
%Degree & Diagrams\\
%\hline
%1 & $D(0,0)$\\
%2 & $D(0,1)$\\
%4 & $D(0,2)$, $D(2,1)$\\
%5 & $D(0,1,1,0)$, $D(2,2)$, $D(1,1,1,1)$, $D(0,3)$, $D(1,2,0)$\\
%7 & $D(5)$, $D(1,4)$, $D(0,4)$\\
%\hline
%\end{tabular}
%\end{center}
\begin{figure}[!ht]
\begin{center}
\begin{tabular}{ccccc}
{\scalebox{.7}{\includegraphics{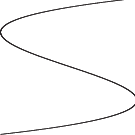}}} &
{\scalebox{.7}{\includegraphics{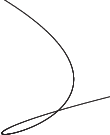}}} &
{\scalebox{.7}{\includegraphics{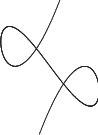}}} &
{\scalebox{.7}{\includegraphics{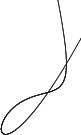}}}&
{\scalebox{.7}{\includegraphics{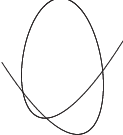}}}\\
$D(0,0)$ & $D(0,1)$ & $D(0,1,1,0)$ & $D(0,2)$ &$D(2,1)$
\end{tabular}
\begin{tabular}{cccc}
{\scalebox{.7}{\includegraphics{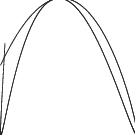}}}&
{\scalebox{.7}{\includegraphics{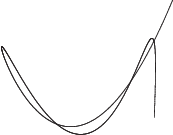}}}&
{\scalebox{.7}{\includegraphics{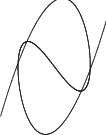}}}&
{\scalebox{.7}{\includegraphics{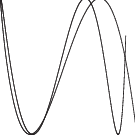}}}
\\
$D(1,2,0)$&$D(0,3)$ &$D(2,2)$ & $D(0,4)$
\end{tabular}
\end{center}
\caption{\small Algebraic degree of a few plane diagrams}\label{fig:ddd}
\end{figure}
\end{lemma}
\Pf
These plane diagrams are obtained with the following parametrisations --
here we use the monic Chebyshev polynomials (also called Dickson polynomials) defined by $T_n(2\cos x) = 2 \cos nx$:   
  \bi
  \item $D(0,0)$: $(T_3,T_1)$
  \item $D(0,1)$ : $(T_3, T_2 - \frac 32 T_1)$
  \item $D(0,1,1,0)$: $(T_3,t^5-4t^3+4t)$
  \item $D(0,2)$:  $(T_3,T_4+3T_2-4T_1)$
\item $D(2,1)$: $(T_3, T_4 - T_2 + \frac 14 T_1)$
 \item $D(0,3)$: $(T_3, t^5-\frac 94 t^4- t^3+\frac{13}4 t^2+\frac 18 t)$
 \item $D(1,2,0)$: $(T_3, T_5(6t/5 + 1/2))$
 \item $D(2,2)$: $(T_3, T_5-\frac{13}{12}T_1)$
 \item $D(0,4)$: $(T_3,T_7(-\frac 32 t + 1))$
\item $D(5)$: $(T_3,P_7)$ where
$P_7 = {t}^{7}-{\frac {93659}{10000}\,{t}^{5}}-{\frac {13549}{5000}\,{t}^{4}}
+{\frac {16453}{1000} \,{t}^{2}}+{\frac {57281}{1000}\,t}$
\item $D(1,4)$: $(T_3,Q_7)$ where
$Q_7 = {t}^{7}-{\frac {84497}{10359}}\, t^5-{\frac {47123}{18875}}\,t^4 +
{\frac {54585}{2759}}\,t^3+{\frac {85741}{7122}}\,t^2-
{\frac {208133}{17097}}\,t -{\frac {242151}{26615}}$
\end{itemize}
%Corollary  \ref{cor:lexdeg1}
It is shown in \cite{BKP2} that the degree is minimal for $D(n,m), \, n,m \geq 0$.
In the case of $D(0,1,1,0)$, every line passing through the two crossing
points meets the curve at 5 points at least, and therefore the degree is at least 5, which is the degree of our parametrisation.
\EPf
\begin{remark}\label{exa:d0n}
One can prove
using {\em dessins d'enfants} (see for example \cite{O4})
that the algebraic
degree of the plane diagram $D(0,n)$ is precisely
$\pent{3n}2+1$.
\end{remark}

\subsection{Two-bridge knots with crossing number at most 6}
\begin{proposition}
  The lexicographic degrees of all two-bridge knots with crossing
number at most $6$ are given in
  Table \ref{tab6}.
\begin{table}[!ht]
\begin{center}
\begin{tabular}{ccrrrrr}
Name & Fraction & Conway Not. & Lex. deg. & Cheb. deg. & diagram & Constr. \\
\hline
$3_{1}$   & $3$ & $ C(3)$ & $(3,4,5)$ & $4$ & $C(3)$ & $D(3)$\\
$4_{1}$  & $5/2$ & $ C(2, 2)$ & $(3,5,7)$ & $5$ & $C(2,2)$ & $D(2,2)$\\
$5_{1}$  & $5$ & $ C(5)$ & $(3,7,8)$ & $7$ & $C(5)$ & $D(5)$\\
$5_{2}$  & $7/2$ & $ C(3, 2)$ & $(3,7,8)$ & $7$ & $C(3,1,1)$ & $D(2,0)+\TT$\\
$6_{1}$  & $9/2$ & $ C(4, 2)$ & $(3,8,10)$ & $8$ & $C(4,2)$ & $D(3, 0)+\TT$\\
$6_{2}$  & $11/3$ & $ C(3, 1, 2)$ & $(3,7,11)$ & $8$ & $C(3,1,2)$ & $D(2,1)+\TT$\\
$6_{3}$  & $13/5$ & $ C(2, 1, 1, 2)$ & $(3,7,11)$ & $7$ & $C(2,1,1,2)$ & $D(0,0)+2\TT$\\
\hline
\end{tabular}
\caption{\small  Lexicographic degree of two-bridge knots with crossing number at most 6}
\label{tab6}
\end{center}
\end{table}
\end{proposition}
\Pf
The knots  $3_1$, and $5_1$ are torus knots, and the knots
$4_1$, $5_2$, and $6_1$ are twist knots. Hence
their lexicographic degrees are
computed in
\cite{BKP1}.
The knots $6_2$ and $6_3$ admit parametrisations with  $b=N+1$, hence
their lexicographic degree is $(3,7,11)$.
\EPf
\subsection{Two-bridge knots with crossing number 7}
\begin{proposition}
    The lexicographic degrees of all two-bridge knots with crossing
number $7$ are given in
  Table \ref{tab7}.
\begin{table}[!ht]
\begin{center}
\begin{tabular}{ccrrrrr}
Name & Fraction & Conway Not. & Lex. deg. & Cheb. deg. & diagram & Constr. \\
\hline
$7_{1}$  & $7$ & $ C(7)$ & $(3,10,11)$ & $10$ & $C(7)$ & $D(7)$\\
$7_{2}$  & $11/2$ & $ C(5, 2)$ & $(3,10,11)$ & $10$ & &Cheb.\\
$7_{3}$  & $13/3$ & $ C(4, 3)$ & $(3,10,11)$ & $10$ & & Cheb.\\
$7_{4}$  & $15/4$ & $ C(3, 1, 3)$ & $(3,8,13)$ & $10$ & $C(3,1,3)$ & $D(1)+2 \TT$\\
$7_{5}$  & $17/5$ & $ C(3, 2, 2)$ & $(3,10,11)$ & $10$ & $C(2, 1, 1,-4)$ & $D(5)+\TT$\\
$7_{6}$  & $19/7$ & $ C(2, 1, 2, 2)$ & $(3,8,13)$ & $10$ & & $D(1)+2\TT$\\
$7_{7}$  & $21/8$ & $ C(2, 1, 1, 1, 2)$ & $(3,8,13)$ & $8$ && Cheb.\\
\hline
\end{tabular}
\caption{\small  Lexicographic degrees of two-bridge knots with crossing number 7}
\label{tab7}
\end{center}
\end{table}
\end{proposition}
\Pf
The lexicographic degree of such a knot  is $(3,8,13)$ or $(3,10,11)$.
The torus knot $7_1$ and the twist knots $7_2$ and $7_3$ have lexicographic
degree $(3,10,11)$, see \cite{BKP2}. The Fibonacci knot $7_7$ has degree
$(3,8,13)$, see \cite{KP4}. The knots $7_4$ and $7_6$ are obtained from
$C(1)$ by $\TT$-augmentation. There degrees is $(3,8,13)$.
The alternating diagram of the knot $7_5$ is $ C(3,2,2).$
By Proposition \ref{prop:Bezout} the
the degree of this diagram is at least $(3,10, 11)$.
Since a non-alternating diagrams of $7_5$ has at least $8$ crossings,
we see that its degree is at least $(3,10,11)$. Hence
the lexicographic degree of $7_5$ is at least $(3,10,11)$.
\EPf

\subsection{Two-bridge knots with crossing number 8}
\begin{proposition}
      The lexicographic degrees of all two-bridge knots with crossing
number  $8$ are given in
  Table \ref{tab8}.
\begin{table}[!ht]
\begin{center}
\begin{tabular}{ccrrrrr}
Name & Fraction & Conway Not. & Lex. deg. & Cheb. deg. & diagram & Constr. \\
\hline
$8_{1}$  & $13/2$ & $ C(6, 2)$ & $(3,11,13)$ & $11$ & & Cheb.\\
$8_{2}$  & $17/3$ & $ C(5, 1, 2)$ & $(3,10,14)$ & $11$&  & $D(4,1)+\TT$\\
$8_{3}$  & $17/4$ & $ C(4, 4)$ & $(3,11,13)$ & $11$ & & Cheb.\\
$8_{4}$  & $19/4$ & $ C(4, 1, 3)$ & $(3,10,14)$ & $11$ &$C(4,1,2,1)$ & $D(2, 0)+2\TT$\\
$8_{6}$  & $23/7$ & $ C(3, 3, 2)$ & $(3,10,14)$ & $11$ & $C(2,2,1,-4)$ & $D(1,2)+2\TT$\\
$8_{7}$  & $23/5$ & $ C(4, 1, 1, 2)$ & $(3,10,14)$ & $10$ && Cheb.\\
$8_{8}$  & $25/9$ & $ C(2, 1, 3, 2)$ & $(3,10,14)$ & $10$ && Cheb.\\
$8_{9}$  & $25/7$ & $ C(3, 1, 1, 3)$ & $(3,10,14)$ & $11$ & & $D(5)+\TT$\\
$8_{11}$  & $27/8$ & $ C(3, 2, 1, 2)$ & $(3,10,14)$ & $11$ &  & $D(2, 0)+2\TT$\\
$8_{12}$  & $29/12$ & $ C(2, 2, 2, 2)$ & $(3,11,13)$ & $11$ & & Cheb.\\
$8_{13}$  & $29/8$ & $ C(3, 1, 1, 1, 2)$ & $(3,10,14)$ & $10$ & & Cheb.\\
$8_{14}$  & $31/12$ & $ C(2, 1, 1, 2, 2)$ & $(3,10,14)$ & $11$ & &$D(2, 0)+2\TT$\\
\hline
\end{tabular}
\caption{\small  Lexicographic degrees of two-bridge knots with crossing number 8}
\label{tab8}
\end{center}
\end{table}
\end{proposition}
\Pf
The lexicographic degree of such a knot  is $(3,10,14)$ or $(3,11,13)$.
The lexicographic degree $(3,11,13)$ of the twist knots $8_1$ and $8_3$
has been obtained in \cite{BKP2}.
Combining Propositions \ref{prop:Bezout} with Chebyshev knots we
obtain the following.
\bi
\item
The knots $8_7$, $8_8$, and $8_{13}$ have minimal lexicographic degree
$(3,10,14)$, obtained as Chebyshev knots.
\item
The plane projection of $8_2=C(5,1,2)$ reduces
to $D(4,1)$ by  $\TT$-reduction.
Since $D(4,1)$ has algebraic degree $7$,
the diagram $D(5,1,2)$ has algebraic degree $10$.
Consequently, $8_2$ has  lexicographic degree $(3,10,14)$.
\item The plane projection of $8_9=C(3,1,1,3)$
  reduces to
  $D(5)$ by
$\TT$-reduction. Hence
  the algebraic degree of $D(3,1,1,3)$ is $10$, and
 $8_9$ has  lexicographic degree $(3,10,14)$.
\item
$D(2,0)$ is obtained by two successive $\TT$-reductions
from the plane projections of diagrams of  $8_4$, $8_{11}$ and $8_{14}$.
Consequently,  $8_4$, $8_{11}$ and $8_{14}$ have lexicographic degree $(3,10,14)$.
%\item The plane diagram
%$D(0,2)$ is obtained by two successive $\TT$-reductions from the
%  plane projections of the alternating diagrams of
%  $8_{14}$.
%  Consequently,  $8_{14}$ has lexicographic degree $(3,10,14)$.
\item Using two $\TT$-reductions, the plane diagram
$D(2,2,1,4)$ reduces to $D(1,2)$, which has algebraic degree 4.
By Proposition \ref{prop:perturb}, the plane diagram
$D(2,2,1,4)$ has algebraic degree 10, and the knot $8_6$ has
lexicographic degree $(3,10,14)$.
\item
The knot $8_{12}$ admits only three simple diagrams with 9 crossings or fewer:
$C(2,2,2,2)$, $C(2, 1, 1, -3, -2)$ and $C(2,2,1,1,-3)$.
By Proposition \ref{prop:Bezout}, the plane diagram
$D(2,2,2,2)$ has
degree at least 11. The plane diagrams
$D(2, 1, 1, 3, 2)$ and $D(2,2,1,1,3)$ reduce,
with two $\TT$-reductions, to $D(3,0)$ or $D(0,3)$ that have
pseudoholomorphic degree  $5$.
By Proposition \ref{prop:perturb},
the lexicographic degree of $8_{12}$ is then $(3,11,13)$.
\EPf
\end{itemize}

The next result shows that the knot $8_6$ is the first example of a knot for which the
lexicographic degree cannot be obtained for the alternating diagram.
\begin{proposition}\label{prop:inters2}
Let $t\mapsto  (P(t), Q(t))$,  be a parametrisation of the diagram
$D(2,3,3)$, where  $\deg P=3$. Then $\deg Q\ge 11$.
\end{proposition}
\Pf
Without loss of generality, we may assume that $P(t)$ is positive
for $t$ large enough, and $\deg Q \not \equiv 0 \Mod 3$.
Let us denote by $C$ the complex algebraic curve image of the map
$t\in \CC \mapsto  (P(t), Q(t))\in\CC^2$.
The curve $C$ has exactly $\deg Q -1$ nodes in $\CC^2$ and then $\deg Q \geq 10$.
Let us suppose that $\deg Q=10$.
Since $C$ has 8 real crossings, it also has a ninth solitary real point.
\begin{figure}[!ht]
\begin{center}
{\scalebox{.6}{\includegraphics{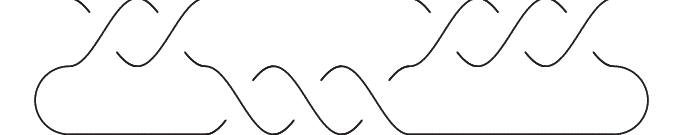}}}
\end{center}
\caption{\small  $C(2,3,3)$}
\end{figure}
We see that there are exactly eight possibilities for the
$\LL$-scheme realised by $C$ (here we use the notations of
\cite[Section 2.2]{BKP2}):
\begin{eqnarray*}
  \supset_2 & \subset_1 \times_2 \times_2 \times_1 \times_1
    \times_1 \times_2 \times_2 \times_2 \supset_1 \bullet_1 &
    \subset_1 \supset_1 \subset_1\\
  \supset_1 & \subset_1 \times_2 \times_2 \times_1 \times_1
    \times_1 \times_2 \times_2 \times_2 \supset_1 \bullet_1 &
    \subset_2 \supset_2 \subset_2\\
   \supset_2 & \bullet_1 \subset_1 \times_2 \times_2 \times_1
     \times_1 \times_1 \times_2 \times_2 \times_2 \supset_1 &
     \subset_1 \supset_1 \subset_1\\
   \supset_1 & \bullet_1 \subset_1 \times_2 \times_2 \times_1
     \times_1 \times_1 \times_2 \times_2 \times_2 \supset_1 &
     \subset_2 \supset_2 \subset_2\\
  \supset_2 & \subset_1 \times_2 \times_2 \times_1 \times_1
    \times_1 \times_2 \times_2 \times_2 \supset_1 \bullet_2 &
    \subset_1 \supset_1 \subset_1\\
  \supset_1 & \subset_1 \times_2 \times_2 \times_1 \times_1
    \times_1 \times_2 \times_2 \times_2 \supset_1 \bullet_2 &
    \subset_2 \supset_2 \subset_2\\
   \supset_2 & \bullet_2 \subset_1 \times_2 \times_2 \times_1
     \times_1 \times_1 \times_2 \times_2 \times_2 \supset_1 &
     \subset_1 \supset_1 \subset_1\\
   \supset_1 & \bullet_2 \subset_1 \times_2 \times_2 \times_1
     \times_1 \times_1 \times_2 \times_2 \times_2 \supset_1 &
     \subset_2 \supset_2 \subset_2
\end{eqnarray*}
We compute all corresponding braids and obtain
\begin{eqnarray*}
&&\bb_1=  \sigma_2^{-1} \sigma_1^{-1} \sigma_2^{-1} \sigma_1^{-3}
    \sigma_2^{-3} \sigma_1^{-3} (\sigma_1 \sigma_2 \sigma_1)^4, \\
&&\bb_2= \sigma_1^{-1} \sigma_2^{-2} \sigma_1^{-3} \sigma_2^{-3}
   \sigma_1^{-2} \sigma_2^{-1} \sigma_1 \sigma_2^{-1} (\sigma_1
   \sigma_2 \sigma_1)^4, \\
&&\bb_3= \sigma_2^{-1} \sigma_1^{-1} \sigma_2 \sigma_1^{-1} \sigma_2^{-2\
  } \sigma_1^{-3} \sigma_2^{-3} \sigma_1^{-2} (\sigma_1
   \sigma_2 \sigma_1)^4, \\
&&\bb_4=\sigma_1^{-2} \sigma_2^{-2} \sigma_1^{-3} \sigma_2^{-3}
   \sigma_1^{-1} \sigma_2^{-1} \sigma_1 \sigma_2^{-1} (\sigma_1
   \sigma_2 \sigma_1)^4, \\
%\end{eqnarray*}% remettre si formule trop longue
%\begin{eqnarray*}
&&\bb_5= \sigma_2^{-1} \sigma_1^{-1} \sigma_2^{-1} \sigma_1^{-3}
   \sigma_2^{-3} \sigma_1^{-1} \sigma_2^{-1} \sigma_1
   \sigma_2^{-1} \sigma_1^{-1} \sigma_2 \sigma_1^{-1} (\sigma_1
   \sigma_2 \sigma_1)^4, \\
&&\bb_6= \sigma_1^{-1} \sigma_2^{-2} \sigma_1^{-3} \sigma_2^{-3}
   \sigma_1^{-1} \sigma_2^{-1} \sigma_1 \sigma_2^{-2} (\sigma_1
   \sigma_2 \sigma_1)^4, \\
&&\bb_7=  \sigma_2^{-2} \sigma_1^{-1} \sigma_2^{-1} \sigma_1^{-3}
    \sigma_2^{-3} \sigma_1^{-2} (\sigma_1 \sigma_2 \sigma_1)^4, \\
&&\bb_8= \sigma_1^{-1} \sigma_2^{-1} \sigma_1 \sigma_2^{-1} \sigma_1^{-1\
  } \sigma_2^{-1} \sigma_1^{-3} \sigma_2^{-3} \sigma_1^{-1}
   \sigma_2^{-1} \sigma_1 \sigma_2^{-1} (\sigma_1 \sigma_2
   \sigma_1)^4.
\end{eqnarray*}
These 8 braids have integer length 0, and none of them is the trivial
braid. Hence the result follows from Proposition \ref{prop:braid}.
\EPf
\pn
Combining Propositions \ref{prop:Bezout} and \ref{prop:inters2},
we obtain
\begin{corollary}
The lexicographic degree of $8_6 = C(2,3,3)$ is not obtained for the alternating diagram.
\end{corollary}
This phenomenon will appear with other knots (see Table \ref{table:noalt}).

\subsection{Two-bridge knots with crossing number 9}
\begin{proposition}
        The lexicographic degrees of all two-bridge knots with crossing
number  $9$ are given in
  Table \ref{table:9}.
\begin{table}[!ht]
\begin{center}
\begin{tabular}{ccrrrrr}
Name & Fraction & Conway Not. & Lex. deg. & Cheb. deg. & diagram & Constr. \\
\hline
 $9_{1}$  & $9$ & $ C(9)$ & $(3,13,14)$ & $13$ && Cheb.\\
 $9_{2}$  & $15/2$ & $ C(7, 2)$ & $(3,13,14)$ & $13$ && Cheb.\\
 $9_{3}$  & $19/3$ & $ C(6, 3)$ & $(3,13,14)$ & $13$ && Cheb.\\
 $9_{4}$  & $21/4$ & $ C(5, 4)$ & $(3,13,14)$ & $13$ && Cheb.\\
 $9_{5}$  & $23/4$ & $ C(5, 1, 3)$ & $(3,11,16)$ & $13$&$C(5,1,2,1)$ & $D(3, 0)+2\TT$\\
 $9_{6}$  & $27/5$ & $ C(5, 2, 2)$ & $(3,13,14)$ & $13$ & &Cheb.\\
 $9_{7}$  & $29/9$ & $ C(3, 4, 2)$ & $(3,13,14)$ & $13$ & &Cheb.\\
 $9_{8}$  & $31/11$ & $ C(2, 1, 4, 2)$ & $(3,11,16)$ & $13$ &$C(2,1,4,1,1)$&  $D(1,2,0)+2\TT$\\
 $9_{9}$  & $31/7$ & $ C(4, 2, 3)$ & $(3,13,14)$ & $13$ & &Cheb.\\
 $9_{10}$  & $33/10$ & $ C(3, 3, 3)$ & $(3,11,16)$ & $13$ &$C(3,2,1,-4)$& $D(0,1)+3\TT$\\
 $9_{11}$  & $33/7$ & $ C(4,1,2,2)$ & $(3,10,17)$ & $13$ && $D(3)+2\TT$\\
 $9_{12}$  & $35/8$ & $ C(4, 2,1, 2)$ & $(3,11,16)$ & $13$& & $D(3, 0)+2\TT$\\
 $9_{13}$  & $37/10$ & $ C(3, 1, 2, 3)$ & $(3,10,17)$ & $13$ & &$D(1,2)+2\TT$\\
 $9_{14}$  & $37/8$ & $ C(4, 1, 1, 1, 2)$ & $(3,11,16)$ & $11$ & & $D(3,0)+2\TT$\\
 $9_{15}$  & $39/16$ & $ C(2, 2, 3, 2)$ & $(3,11,16)$ & $13$ & $C(2,2,2,1,-3)$&
$D(1,0)+3\TT$\\
 $9_{17}$  & $39/14$ & $ C(2, 1, 3, 1, 2)$ & $(3,10,17)$ & $11$ & &$D(3)+2\TT$\\
 $9_{18}$  & $41/12$ & $ C(3,2, 2, 2)$ & $(3,13,14)$ & $13$ && Cheb.\\
 $9_{19}$  & $41/16$ & $ C(2, 1, 1, 3, 2)$ & $(3,11,16)$ & $11$ & &$D(3, 0)+2\TT$\\
 $9_{20}$  & $41/11$ & $ C(3, 1, 2, 1, 2)$ & $(3,10,17)$ & $13$ & &$D(3)+2\TT$\\
 $9_{21}$  & $43/12$ & $ C(3,1,1,2,2)$ & $(3,11,16)$ & $13$ & &$D(3, 0)+2\TT$\\
 $9_{23}$  & $45/19$ & $ C(2, 2, 1, 2, 2)$ & $(3,10,17)$ & $13$ & &$D(0, 0)+3\TT$\\
 $9_{26}$  & $47/13$ & $ C(3, 1, 1, 1, 1, 2)$ & $(3,11,16)$ & $11$ & &$D(3)+2\TT$\\
 $9_{27}$  & $49/18$ & $ C(2, 1, 2, 1, 1, 2)$ & $(3,10,17)$ & $13$ & &$D(3)+2\TT$\\
$9_{31}$   & $55/21$ & $ C(2, 1, 1, 1, 1, 1, 2)$ & $(3,10,17)$ & $10$ && Cheb.\\
\hline
\end{tabular}
\caption{\small  Lexicographic degree of two-bridge knots with crossing number 9}
\label{table:9}
\end{center}
\end{table}
\end{proposition}
\Pf
The lexicographic degree
of such a knot is
$(3,10,17), (3,11,16)$, or $(3,13,14)$.
Furthermore, any  diagram with at least 11 crossings
has degree $(3,13,14)$ at least.
It is proved in  \cite{KP4}
that $9_{31}$ is the harmonic Fibonacci
knot $(T_3,T_{10},T_{17})$.
The torus knot $9_1$ and the twist knots $9_2, 9_3, 9_4,$
have lexicographic degree $(3,13,14)$, see \cite{BKP2}.
For the remaining knots, we proceed as follows.
\bi
\item By $\TT$-reduction, the diagram $D(2,2,1,2,2)$
reduces  to   $D(0,0)$, that   has algebraic degree 1.
  We deduce that the knot $9_{23} = C(2,2,1,2,2)$
has lexicographic degree $(3,10,17)$.
\item The alternating diagrams of $9_{11}, 9_{13}$, $9_{17}$, $9_{20}$, $9_{26}$
and $9_{27}$ can be reduced to
$D(3)$ by two $\TT$-reductions. Their lexicographic degree is then $(3,10,17)$.
\item The plane alternating diagram of $9_8$ is reduced to
$D(1,3,2)$ by $\TT$-reduction. The algebraic degree of $D(1,3,2)$ is at most
the degree of $D(4,2)$, that is 8. On the other hand, the plane projection of the diagram $C(2,1,4,1,1)$ can be reduced to $D(1,2,0)$ that has degree 8.
\item
The plane alternating diagrams of the knots $9_5$, $9_{12}$, $9_{14}$, $9_{19}$ and
$9_{21}$ can be reduced by two $\TT$-reductions to $D(3,0)$. Hence these
diagrams have algebraic degree 11.
On the other hand,
any other diagram of these knots will be non-alternating with at least
10 crossing points.
Hence the lexicographic degree of these knots is then $(3,11,16)$.
\item The alternating diagram of $9_{15}$ is $C(2,2,3,2)$.
From Proposition \ref{prop:Bezout}, its lexicographic degree is at least
$(3,13,14)$.
Any other non alternating diagram of $9_{15}$ will have 10 or more crossings.
Consider the diagram $C(2,2,2,1,-3)$ of $9_{15}$. Its projection $D(2,2,2,1,3)$
can be reduced to $D(1,0)$ by three $\TT$-reductions. Consequently $9_{15}$ has degree $(3,11,16)$.
\item  The alternating diagram of $9_6$ is $C(5,2,2)$.
From Proposition \ref{prop:Bezout}, its lexicographic degree is at least
$(3,13,14)$. The only diagrams of $9_6$ having 10 crossings are
  $C(2,1,1,-6)$ and $C(5,1,1,-3)$, whose plane diagrams reduce to $D(7)$
by $\TT$-reductions.
Hence the lexicographic degree of $9_6$ is  $(3,11,16)$.
\item  The alternating diagram of $9_7$ is $C(3,4,2)$.
From Proposition \ref{prop:Bezout}, its lexicographic degree is at least
$(3,13,14)$. The only diagrams of $9_7$ having 10 crossings are
$C(2,3,1,-4)$, and $C(3,3,1,-3)$.
The  plane diagrams $D(2,3,1,4)$ reduces to $D(2,2,3)$ and
$D(3,3,1,3)$ to $D(3,2,2)$ by a $\TT$-reduction.
Their degrees are at least $14$ by Proposition \ref{prop:inters2}.
\item  The alternating diagram of $9_9$ is $C(4,2,3)$.
From Proposition \ref{prop:Bezout}, its lexicographic degree is at least
$(3,13,14)$. The only diagrams of $9_9$ having 10 crossings are
 $C(3,1,1,-5)$ and $C(4,1,1,-4)$,  whose plane diagrams reduce to $D(7)$
by $\TT$-reductions. Their lexicographic degrees are then $(3,13,14)$.
\item The alternating diagram of $9_{18}$ is $C(3,2,2,2)$.
From Proposition \ref{prop:Bezout}, its lexicographic degree is at least
$(3,13,14)$. The only diagrams of $9_{18}$ having 10 crossings are
$C(3, 1, 1, -3, -2)$ whose plane projection reduces to $D(5,2)$,
$C(2, 2, 1, 1, -4)$ whose plane projection reduces to $D(2,5)$,
$C(2, 1, 1, -3, -3)$ whose plane projection reduces to $D(4,3)$,
and $C(3,2,1,1,-3)$, whose plane projection reduces to $D(3,4)$.
By Proposition \ref{prop:even}, the degree of these four plane
diagrams with seven crossings is at least $10$, so the degree of the
four plane diagrams with 10 crossings is at least 13 by Proposition
\ref{prop:flip}.
\item  The alternating diagram of $9_{10}$ is $C(3,3,3)$. Suppose that
  there exists  a polynomial parametrisation  $\gamma:t\mapsto  (P(t), Q(t))$
  of the plane diagram $D(3,3,3)$ with $\deg(P)=3$ and
  $\deg(Q)=10$. We denote by $C=\gamma(\CC)$.
Since the curve $C$ has 9 real crossings,
  it has no additional nodes.
The braid associated to $C$ is
$$b_C=\sigma_1^{-1}\sigma_2^{-1}\sigma_1^{-2}\sigma_2^{-3}\sigma_1^{-3}\sigma_2^{-2}(\sigma_1\sigma_2\sigma_1)^4.$$
Since this braid is not the trivial braid, we obtain a
contradiction. Hence the alternating diagram $C(3,3,3)$ has degree at
least $(3,11,16)$. On the other hand, the projection of the
diagram $C(3,2,1,-4)$ of $9_{10}$ reduces
to $D(2,2)$. Since this latter has algebraic degree 5,
we deduce that $9_{10}$   has lexicographic degree $(3,11,16)$.
\EPf
\end{itemize}
 \begin{figure}[!ht]
 \begin{center}
 \begin{tabular}{ccc}
 {\scalebox{.45}{\includegraphics{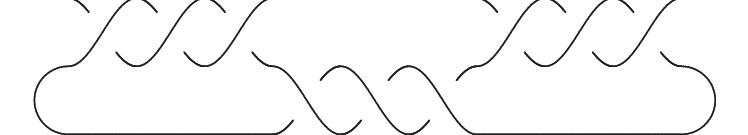}}}&
 {\scalebox{.45}{\includegraphics{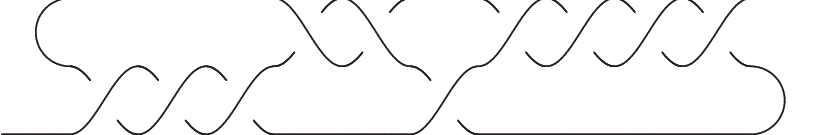}}}\\
 $C(3,3,3)$&$C(3,2,1,-4)$.
 \end{tabular}
 \vspace{-10pt}
 \caption{\small Two diagrams of $9_{10}$}
 \end{center}
 \end{figure}
\subsection{Two-bridge knots with crossing number 10 or 11}\label{sec:10-11}
The lexicographic degrees
of  the torus knot $C(11)$ and the twist knots $C(8,2)$, $C(9,2)$, $C(8,3)$,
$C(6,4)$, $C(7,4)$ and $C(6,5)$ have been established in  \cite{BKP2}.
For the 129
remaining knots with 10 or 11 crossings, we
simply sketch all computations.
For only 11 knots among the 186 knots with 11  crossings or fewer ---
 $ 10_{ 11}$, $ 10_{ 13}$, $11a_{ 98}$, $11a_{166}$, $11a_{230}$, $11a_{235}$, $11a_{238}$, $11a_{311}$, $11a_{335}$, $11a_{359}$ and $11a_{365}$ ---
the lower bounds differ from the upper bounds in the strategy described
in Section \ref{sec:general}, i.e. one has to go through step 4.
The projections of all the  corresponding diagrams reduce by $\TT$-reduction
to a finite list of eleven plane diagrams:
\bi
\item
$D(3,3,3)$ and $D(3,3,4)$, that have degree $13$ at least,
\item
$D(3,3,5)$ and $(3,5,3)$, that have degree $14$ at least,
\item
$D(3,3,6)$, $D(3,5,4)$, $D(3,2,3,4)$, $D(3,2,3,5)$, $D(3,2,5,3)$, that
have degree $16$ at least,
\item
$D(3, 3, 2, 5)$ and $D(4, 2, 3, 4)$ that have degree $(3,17)$ at least.
\ei
These results have been obtained by computing all possible braids associated
  to hypothetical plane curves of degree $b<b_0$ that are
  $\LL$-isotopic to the diagram, and checking,
like in Proposition \ref{prop:inters2}, if these braids satisfy Proposition \ref{prop:braid}.

\section{Conclusion}
We list in Table \ref{table:all} the lexicographic degrees of the first 186
two-bridge knots. We only write $b$, bearing in mind that the
corresponding lexicographic degree is $(3,b,3N-b)$.
Details of our results will be available in
\href{https://webusers.imj-prg.fr/~pierre-vincent.koseleff/knots/2bk-lexdeg.html}
{https://.../2bk-lexdeg.html}
\begin{table}[!ht]
\begin{center}
\begin{tabular}{||rc|rc|rc|rc|rc|rc||}
\hline
Name& Deg.&
Name& Deg.&
Name& Deg.&
Name& Deg.&
Name& Deg.&
Name& Deg.\\
\hline
$  3_{  1}$ &4&$  4_{  1}$ &5&$  5_{  1}$ &7&$  5_{  2}$ &7&$  6_{  1}$ &8&$  6_{  2}$ &7\\
$  6_{  3}$ &7&$  7_{  1}$ &10&$  7_{  2}$ &10&$  7_{  3}$ &10&$  7_{  4}$ &8&$  7_{  5}$ &10\\
$  7_{  6}$ &8&$  7_{  7}$ &8&$  8_{  1}$ &11&$  8_{  2}$ &10&$  8_{  3}$ &11&$  8_{  4}$ &10\\
$  8_{  6}$ &10&$  8_{  7}$ &10&$  8_{  8}$ &10&$  8_{  9}$ &10&$  8_{ 11}$ &10&$  8_{ 12}$ &11\\
$  8_{ 13}$ &10&$  8_{ 14}$ &10&$  9_{  1}$ &13&$  9_{  2}$ &13&$  9_{  3}$ &13&$  9_{  4}$ &13\\
$  9_{  5}$ &11&$  9_{  6}$ &13&$  9_{  7}$ &13&$  9_{  8}$ &11&$  9_{  9}$ &13&$  9_{ 10}$ &11\\
$  9_{ 11}$ &10&$  9_{ 12}$ &11&$  9_{ 13}$ &10&$  9_{ 14}$ &11&$  9_{ 15}$ &11&$  9_{ 17}$ &10\\
$  9_{ 18}$ &13&$  9_{ 19}$ &11&$  9_{ 20}$ &10&$  9_{ 21}$ &11&$  9_{ 23}$ &10&$  9_{ 26}$ &10\\
$  9_{ 27}$ &10&$  9_{ 31}$ &10&$ 10_{  1}$ &14&$ 10_{  2}$ &13&$ 10_{  3}$ &14&$ 10_{  4}$ &13\\
$ 10_{  5}$ &13&$ 10_{  6}$ &13&$ 10_{  7}$ &13&$ 10_{  8}$ &13&$ 10_{  9}$ &13&$ 10_{ 10}$ &13\\
$ 10_{ 11}$ &13&$ 10_{ 12}$ &13&$ 10_{ 13}$ &14&$ 10_{ 14}$ &13&$ 10_{ 15}$ &13&$ 10_{ 16}$ &11\\
$ 10_{ 17}$ &13&$ 10_{ 18}$ &13&$ 10_{ 19}$ &13&$ 10_{ 20}$ &13&$ 10_{ 21}$ &13&$ 10_{ 22}$ &13\\
$ 10_{ 23}$ &13&$ 10_{ 24}$ &13&$ 10_{ 25}$ &13&$ 10_{ 26}$ &13&$ 10_{ 27}$ &13&$ 10_{ 28}$ &11\\
$ 10_{ 29}$ &11&$ 10_{ 30}$ &11&$ 10_{ 31}$ &13&$ 10_{ 32}$ &13&$ 10_{ 33}$ &11&$ 10_{ 34}$ &13\\
$ 10_{ 35}$ &14&$ 10_{ 36}$ &13&$ 10_{ 37}$ &13&$ 10_{ 38}$ &11&$ 10_{ 39}$ &13&$ 10_{ 40}$ &13\\
$ 10_{ 41}$ &11&$ 10_{ 42}$ &11&$ 10_{ 43}$ &11&$ 10_{ 44}$ &11&$ 10_{ 45}$ &11&$11a_{ 13}$ &14\\
$11a_{ 59}$ &14&$11a_{ 65}$ &14&$11a_{ 75}$ &13&$11a_{ 77}$ &13&$11a_{ 84}$ &13&$11a_{ 85}$ &13\\
$11a_{ 89}$ &13&$11a_{ 90}$ &13&$11a_{ 91}$ &13&$11a_{ 93}$ &13&$11a_{ 95}$ &13&$11a_{ 96}$ &14\\
$11a_{ 98}$ &14&$11a_{110}$ &13&$11a_{111}$ &13&$11a_{117}$ &13&$11a_{119}$ &14&$11a_{120}$ &13\\
$11a_{121}$ &14&$11a_{140}$ &13&$11a_{144}$ &13&$11a_{145}$ &14&$11a_{154}$ &14&$11a_{159}$ &14\\
$11a_{166}$ &14&$11a_{174}$ &13&$11a_{175}$ &13&$11a_{176}$ &13&$11a_{177}$ &13&$11a_{178}$ &13\\
$11a_{179}$ &13&$11a_{180}$ &13&$11a_{182}$ &13&$11a_{183}$ &13&$11a_{184}$ &13&$11a_{185}$ &13\\
$11a_{186}$ &13&$11a_{188}$ &13&$11a_{190}$ &13&$11a_{191}$ &13&$11a_{192}$ &13&$11a_{193}$ &13\\
$11a_{195}$ &14&$11a_{203}$ &13&$11a_{204}$ &13&$11a_{205}$ &13&$11a_{206}$ &13&$11a_{207}$ &13\\
$11a_{208}$ &13&$11a_{210}$ &14&$11a_{211}$ &14&$11a_{220}$ &13&$11a_{224}$ &13&$11a_{225}$ &13\\
$11a_{226}$ &14&$11a_{229}$ &14&$11a_{230}$ &14&$11a_{234}$ &16&$11a_{235}$ &16&$11a_{236}$ &16\\
$11a_{238}$ &16&$11a_{242}$ &16&$11a_{243}$ &16&$11a_{246}$ &16&$11a_{247}$ &16&$11a_{306}$ &13\\
$11a_{307}$ &13&$11a_{308}$ &13&$11a_{309}$ &13&$11a_{310}$ &13&$11a_{311}$ &14&$11a_{333}$ &14\\
$11a_{334}$ &16&$11a_{335}$ &16&$11a_{336}$ &13&$11a_{337}$ &13&$11a_{339}$ &16&$11a_{341}$ &13\\
$11a_{342}$ &16&$11a_{343}$ &14&$11a_{355}$ &16&$11a_{356}$ &13&$11a_{357}$ &13&$11a_{358}$ &16\\
$11a_{359}$ &14&$11a_{360}$ &13&$11a_{363}$ &14&$11a_{364}$ &16&$11a_{365}$ &14&$11a_{367}$ &16\\
\hline
\end{tabular}
\caption{\small Two-bridge knots with crossing number at most 11 and
  their $y$-lexicographic degree}
  %eb
  % 11 or fewer crossings and their $y$-degree}
\label{table:all}
\end{center}
\end{table}
\pn
In Table \ref{table:noalt}, we list all knots for which the algebraic
degrees of their alternating diagrams are greater than their lexicographic
degrees.
The third column of Table
\ref{table:noalt} gives a diagram obtained by a polynomial
parametrisation of lexicographic degree, the fourth column indicates a
construction of the corresponding $xy$-plane diagram (the notation is
explained in  Section
\ref{sec:plane}), the fifth column gives the alternating trigonal diagram of
the knot, and the last column gives a lower bound on its $y$-degree.
\begin{table}[!ht]
\begin{center}
\begin{tabular}{||rrrr|rr||}
\hline
Name & $y$-lex. &  Lex. deg. & Constr. & Alt. diagram & $y$-lex. \\
& degree &  diagram & & & degree $\ge$ \\
\hline
%[[2, 2, 1, -4], [3, 2, 1, -4], [2, 2, 1, -3, -2], [2, 2, 1, -3, -3], [2, 1, 3, 2, 1, -3], [2, 2, 1, -3, -1, -1, -2], [2, 2, 2, 1, -5], [2, 2, 2, 1, -3, -2], [2, 1, 1, 1, 2, 1, -4], [2, 2, 1, -2, -1, -1, -3], [2, 1, 1, -2, -1, -2, -3], [3, 2, 1, -3, -1, -2], [2, 3, 1, -2, -4], [2, 2, 1, -3, -4], [3, 1, 3, 1, -4], [3, 2, 1, -3, -3]]
$8_{  6}$ & 10& $C(2, 2, 1, -4)$& $D(3)+2 \TT$& $C(3, 3, 2)$& 11\\
$9_{ 10}$ & 11& $C(3, 2, 1, -4)$& $D(0, 1)+3 \TT$& $C(3, 3, 3)$& 13\\
$9_{ 15}$& 11& $C(2, 2, 1, -3, -2)$& $D(1, 0)+ 3 \TT$& $C(2, 2, 3, 2)$& 13\\
$10_{ 24}$& 13& $C(2, 2, 1, -3, -3)$& $D(0, 2)+ 3 \TT$& $C(3, 2, 3, 2)$& 14\\
$11a_{ 75}$& 13& $C(2, 1, 3, 2, 1, -3)$&  $D(3)+ 3 \TT$& $C(2, 1, 3, 3, 2)$& 14\\
$11a_{ 84}$& 13& $C(2, 2, 1, -3, -1, -1, -2)$& $D(0, 0)+ 4 \TT$& $C(2, 1, 1, 2, 3, 2)$& 14\\
$11a_{144}$& 13& $C(2, 2, 2, 1, -5)$& $D(3)+3\TT$& $C(4, 3, 2, 2)$& 14\\
$11a_{186}$& 13& $C(2, 2, 2, 1, -3, -2)$&  $D(0, 0)+4\TT$& $C(2, 2, 3, 2, 2)$& 16\\
$11a_{193}$& 13& $C(2, 1, 1, 1, 2, 1, -4)$& $D(3)+3\TT$& $C(3, 3, 1, 1, 1, 2)$& 14\\
$11a_{205}$& 13& $C(2, 2, 1, -2, -1, -1, -3)$& $D(3)+3\TT$& $C(3, 1, 1, 1, 3, 2)$& 14\\
$11a_{208}$& 13& $C(2, 1, 1, -2, -1, -2, -3)$& $D(3)+3\TT$& $C(3, 2, 1, 1, 2, 2)$& 14\\
$11a_{224}$& 13& $C(3, 2, 1, -3, -1, -2)$& $D(0, 0)+4\TT$& $C(3, 3, 2, 1, 2)$& 14\\
$11a_{225}$& 13& $C(2, 3, 1, -2, -4)$& $D(3)+3 \TT$& $C(4, 1, 4, 2)$& 14\\
$11a_{229}$& 14& $C(2, 2, 1, -3, -4)$& $D(0, 3)+ 3\TT$& $C(4, 2, 3, 2)$& 16\\
$11a_{341}$& 13& $C(3, 1, 3, 1, -4)$& $D(3)+3\TT$& $C(3, 1, 4, 3)$& 14\\
$11a_{356}$& 13& $C(3, 2, 1, -3, -3)$& $D(3)+3 \TT$& $C(3, 2, 3, 3)$& 16\\

\hline
\end{tabular}
\caption{\small Knots for which the alternating diagram is not of minimal degree}
\label{table:noalt}
\end{center}
\end{table}

\pn
For $N\geq 12$ and $N+4 \leq b < \pent{3N-1}2$, it could be interesting
to determine the lexicographic degree, as we do not know if $b+c=3N$.
For some knots, it could be interesting to determine explicit
constructions with the lexicographic degree.

%%%%%%%%%%%%%%%%%%%%%%%%%%%%%%%%%%%%%%%%%%%%%%%%%%%%%%%%%%%%%%%%%%%%%%%%%%%%%%
%@

%\input annex-bkp3.tex
%%%%%%%%%%%%%%%%%%%%%%%%%%%%%%%%%%%%%%%%
\pn
\hrule width 7cm height 1pt %depth 3pt
\pn
{\small
Erwan {\sc Brugall\'e}\\
Universit\'e de Nantes,\\
Laboratoire de Math\'ematiques Jean Leray (CNRS LMJL UMR 6629)\\
e-mail: {\tt erwan.brugalle@math.cnrs.fr}
\pn
Pierre-Vincent {\sc Koseleff}\\
Sorbonne Universit\'e (UPMC--Paris 6),\\
Institut de Math\'ematiques de Jussieu (CNRS IMJ-PRG UMR 7586)  \&
Inria-Paris\\
%4, place Jussieu, F-75252 Paris Cedex 05 \\
e-mail: {\tt pierre-vincent.koseleff@upmc.fr}
\pn
Daniel {\sc Pecker}\\
Sorbonne Universit\'e (UPMC--Paris 6),\\
Institut de Math\'ematiques de Jussieu (CNRS IMJ-PRG UMR 7586),\\
%4, place Jussieu, F-75252 Paris Cedex 05 \\
e-mail: {\tt daniel.pecker@upmc.fr}
}
\vfill\hbox{}
\end{document}